\let\ORIlabel\label
\let\ORIrefstepcounter\refstepcounter
\AddToHook{package/hyperref/before}{%
  \let\label\ORIlabel
  \let\refstepcounter\ORIrefstepcounter
}

\documentclass[hidelinks,onefignum,onetabnum]{siamart220329}

\usepackage{amssymb,amsmath,amsopn}
\usepackage{amsfonts}
\usepackage{amstext}
\usepackage{empheq}
\usepackage{epsfig}
\usepackage{color}
\usepackage{graphicx,colordvi}
\usepackage{epstopdf}
\usepackage{float}

\usepackage{geometry}
\usepackage{lipsum}
\usepackage{amsfonts}
\usepackage{graphicx}
\usepackage{epstopdf}
\usepackage{algorithmic}
\usepackage{fancyhdr}

\ifpdf
  \DeclareGraphicsExtensions{.eps,.pdf,.png,.jpg}
\else
  \DeclareGraphicsExtensions{.eps}
\fi

\newsiamremark{remark}{Remark}
\newsiamremark{hypothesis}{Hypothesis}
\crefname{hypothesis}{Hypothesis}{Hypotheses}
\newsiamthm{claim}{Claim}

\headers{Leader-Follower Mean Field LQG Games with Multiplicative Noise}{B-. C. Wang, H. S. Zhang and J-. F. Zhang}

\title{Leader-Follower Mean Field LQG Games with Multiplicative Noise\thanks{Submitted to the editors DATE.
\funding{This work was funded by the National Natural Science Foundation of China under Grants 62433020, 62573266, 62192753 and T2293770, and the Innovative Research Groups of National Natural Science Foundation of China under Grant 61821004.}}}

\author{Bing-Chang Wang \thanks{School of Control Science and Engineering,
Shandong University, Jinan, China
  (\email{bcwang@sdu.edu.cn}).}
\and Huanshui Zhang \thanks{College of Electrical Engineering and Automation,
	Shandong University of Science and Technology, Qingdao, China
	(\email{hszhang@sdu.edu.cn}).}
\and Ji-Feng Zhang \thanks{  Corresponding author.  School of Automation and Electrical Engineering, Zhongyuan University of Technology, Zhengzhou 450007, Henan Province, China; and the State Key Laboratory of Mathematical Sciences, Academy of Mathematics and Systems Science, Chinese Academy of Sciences, Beijing 100190, China (\email{jif@iss.ac.cn}). }
}

\usepackage{amsopn}

\ifpdf
\hypersetup{
	pdftitle={{Leader-Follower} Mean Field LQG Games with Multiplicative Noise},
	pdfauthor={B. Wang}
}
\fi

\begin{document}
	
	\maketitle
	
	\begin{abstract}
		This paper studies open-loop and feedback solutions to leader-follower mean field linear-quadratic-Gaussian games with multiplicative noise by the direct approach. The leader-follower game involves a leader and many followers, where the state and	control weight matrices in their costs are  not limited to be positive definite. From variational analysis with mean field approximations, we obtain a set of  open-loop controls in terms of solutions to mean field forward-backward stochastic differential equations. By applying the matrix maximum principle, a set of decentralized feedback strategies is constructed.  Distinct from traditional works, a cross term has appeared in derivation due to the presence of mean field terms. For open-loop and feedback solutions, the corresponding optimal costs of all players are explicitly given in terms of the solutions to two Riccati equations, respectively.
	\end{abstract}
	
	\begin{keywords}
		Stackelberg game, mean field team, social control, forward-backward stochastic differential equation
	\end{keywords}
	
	\begin{MSCcodes}
		49N80, 91A16, 93E03, 93E20
	\end{MSCcodes}
	
	\section{Introduction}
	\subsection{Background and Motivation}
	Mean field (MF) games have drawn much attention from various disciplines including 
	control theory, applied mathematics and economics  \cite{ll2007}, \cite{bfy2013}, \cite{c2014}, \cite{CD18}.
	In an MF game, the impact of each individual is negligible while the effect of the population is significant.
	The main methodology of MF games is to replace the interactions among agents by population aggregation effect, which structurally models the MF interactions in large population systems. Thus, the high-dimensional multi-agent optimization problem can be transformed into a low-dimensional local optimal control problem for a representative agent \cite{ll2007}, \cite{c2014}. 
	Wide applications have been found in many fields, such as economics \cite{weintraub2008markov}, \cite{bibitem181}, smart grid \cite{TKG20}, engineering  \cite{bibitem7} and social sciences \cite{bibitem5}, \cite{CDL22}.
	As a classical type of MF models,
	mean field linear quadratic Gaussian (MF-LQG) games are intensively studied due to their analytical tractability and close connection to practical applications.
	For works on such kind of problems, readers can refer to \cite{bp2014}, \cite{jm2017}, \cite{hcm2012}, \cite{lz2008}, \cite{tzb2014}, \cite{bj2012a}, \cite{bj2017b}. The pioneering work  \cite{hcm2007a} studied $\epsilon$-Nash equilibrium strategies for MF-LQG games with discounted costs based on the Nash certainty
	equivalence. This approach was then  applied to the cases with long run average costs \cite{lz2008} and with  Markov jump parameters \cite{bj2012a}, respectively.
	For  MF games with  major players, the works \cite{h2010}, \cite{CK17} considered continuous-time LQG games with complete and partial information;  
	\cite{bj2012b} investigated discrete-time LQG games with random parameters;  \cite{BLP13} and \cite{SC2016} focused on the nonlinear case.

	In contrast to the above models, the leader-follower (Stackelberg) game involves a leader-follower structure. Consider a leader-follower game with two layers. One layer of players are defined as leaders with a dominant position and the other players is defined as followers with a subordinate position. The leader has the priority to give a strategy first and then followers seek strategies to minimize their costs with response to the strategies of leaders. According to followers' optimal response, leaders will choose strategies to minimize their costs. 
	Leader-follower games have been widely investigated in the literature (see e.g. \cite{sc1973}, \cite{y02}, \cite{bbs2010}, \cite{XZC15}, \cite{HWX20}).  Recently, leader-follower MF games have attracted great research interest \cite{bcy2015a}, \cite{bj2014a}, \cite{mb2018}, \cite{BAA20}, \cite{YH21}.  
	The work \cite{bcy2015a} considered MF Stackelberg games with delayed instructions. \cite{bj2014a} studied discrete-time hierarchical MF games with
	tracking-type costs and gave the $\varepsilon$-Stackelberg equilibrium. 
	Authors in \cite{mb2018} investigated continuous-time MF-LQG Stackelberg games by the fixed-point method,
	and they asserted that ``complexity brought by coupling among leader and followers makes the use of direct approach almost impossible".
	This work is further generalized to the jump diffusion model \cite{M23}.
	Furthermore, \cite{YH21} 
	investigated feedback strategies of MF Stackelberg games by solving the master equations.

	Distinct from  noncooperative games,  social optimization is a joint decision problem where all players work cooperatively to optimize the social cost. This is a typical class of team decision problem \cite{hch1972}. Authors in \cite{hcm2012} studied social optima in the MF-LQG control, and provided an asymptotic team-optimal solution, which is extended to the case of mixed games in \cite{hn2016}. The work \cite{bj2017b} investigated the MF social optimal problem where the jump parameter appears as a common source of randomness. 
	More investigation can be found in \cite{AM15} for team-optimal control with finite population and partial information, 
	\cite{SNM18} for dynamic 
	collective choice by finding social optima, \cite{SY20} for stochastic dynamic teams and their MF limit, 
	\cite{WHZ21}, \cite{HWY21} for MF teams with uncertainty in drift and volatility, and \cite{NM18} for social control applications in economics.
In addition, see \cite{WLC23} for value-iteration learning in ergodic MF-LQG social control, and \cite{JWS24} for
 online policy iteration in MF Pareto optimal control.

	Normally, there are two routes to solve MF  games and teams. One is called the fixed-point approach \cite{hcm2007a, hcm2012, bfy2013, CD18}, which starts by applying MF approximation and constructing a fixed-point equation. 
	A set of decentralized strategies can be designed by tackling the fixed-point equation together with the optimal response of a representative player. 
	In general, the fixed-point equation is difficult to solve. In addition, when solving the team problem by the fixed-point approach,
	an additional variable (called social impact \cite{hcm2012, bj2017b}) needs to be introduced. This leads to a drastic increase of computational complexity for MF teams with \emph{multiplicative noise} \cite{FHQ21}, \cite{QHX20}.
	Another route is called the direct approach \cite{
		HZ20,ll2007, WZ21}, which takes a path from finite-population to infinite-population systems. 
	By decoupling the  
	Hamiltonian system for $N$-player, one can obtain a centralized strategy which explicitly relies on the state of a player and population state average.
	Applying MF approximations,  the decentralized control can be constructed. By the direct approach, the resulting control is neat and less computation is required, particularly for team problems \cite{WZ21}.
	

	\subsection{Contribution and Novelty}
	This paper considers MF-LQG Stackelberg games with a leader and many followers, where the state and control
	weight matrices in their costs are allowed to  not be positive definite. The leader first give his strategy and then all followers cooperate to optimize the \emph{social cost}, the sum of individual costs. For instance, consider an example of macroeconomic regulation, where the regulator/
government is the leader, and local authorities are followers \cite{P77}.
	The state of the leader appears in both dynamics and cost of each follower. It shows that the dynamics and costs of followers are directly influenced by the behavior of the leader.
	Distinct from  \cite{hn2016} and \cite{mb2018}, our model involves population state average $x^{(N)}$ in both drift and diffusion terms of all players' dynamics,
which imples the leader and followers are fully coupled by the MF term. {Owing to the presence of indefinite cost weights  and multiplicative noise,  the control design and analysis get more difficult.   Convex  analysis is needed for the leader-follower MF-LQG problem.  In particular, the
convex  analysis for leader's problem is challenging, since  the system is driven by a set of coupled
forward-backward stochastic differential equations (FBSDEs). 

 By the terminology of \cite{BO99}, the solutions to Stackelberg
games are mainly divided into open-loop, closed-loop and feedback (closed-loop memoryless)
solutions.
 The 
Stackelberg solution under closed-loop information pattern cannot be solved by utilizing the
standard techniques of optimal control theory (See \cite[p. 376]{BO99}).
However, the feedback solution to Stackelberg LQG
games with strictly convex cost can be determined in the closed form. Compared with the open-loop solution, there exists stronger coupling
among the feedback strategies of the leader and numerous followers in MF games. Additionally, the MF coupling among players bring about more difficulty in strategy design.
 Until now, most previous works focused on open-loop solutions of  MF leader-follower games, and only a few
	works were on feedback and closed-loop solutions. Furthermore, the relationship among different solutions is still unclear.}
	
	In this paper, we study systematically open-loop and feedback solutions to MF  leader-follower games by the direct approach.
	The open-loop solution starts with solving a centralized social control problem for followers, and obtaining a system of high-dimensional  FBSDEs. By MF approximations, a set of open-loop controls of followers is designed in terms of an MF FBSDE.
	After applying followers' strategies, 
	we derive necessary and sufficient conditions for the solvability of the leader's  problem, and then obtain the feedback representation of the open-loop control by decoupling an FBSDE.  
	From perturbation analysis, the proposed strategy is shown to be an $(\varepsilon_1,\varepsilon_2)$-Stackelberg equilibrium. Furthermore, we obtain the optimal costs of players in terms of the solutions to Riccati equations. Next, the feedback solution is investigated for MF Stackelberg games. Distinct from  the open-loop solution, we  presume that the leader has a strategy with the feedback form. With leader's feedback gain fixed, we obtain the feedback strategies of followers by decoupling high-dimensional FBSDEs. Applying the matrix maximum principle with MF
	approximations, we solve the optimal control problem
	for the leader, and then construct a set of decentralized feedback
	strategies for all players.
	By the technique of  completing the square, we show that the proposed decentralized
	strategy is a feedback $(\varepsilon_1,\varepsilon_2)$-Stackelberg equilibrium and give an explicit form of the corresponding costs of players.

	The main contributions of the paper are listed as follows.
	\begin{itemize}
		\item By adopting a direct approach, {we explore the open-loop and feedback solutions to indefinite leader-follower MF games with multiplicative noise. Distinct from  the fixed-point approach}, 
		\emph{no additional terms}  are introduced when  
		MF social control problem is solved for followers.
		
		\item By variational analysis with MF approximations, we obtain an  open-loop asymptotic Stackelberg equilibrium in terms of MF FBSDEs, which can be implemented offline.	
		
		\item  
		By decoupling high-dimensional FBSDEs and applying the matrix maximum principle,  a set of decentralized feedback
		strategies is constructed. Distinct from traditional works,
		 a cross term has appeared in deriving  feedback
		strategies due to the presence of MF coupling.
		
		
	\end{itemize}
	
	%
	
	\subsection{Organization and Notation}
	
	The paper is organized as follows. In Section 2, we formulate the problem of MF-LQG leader-follower games with  multiplicative noise. 
	In Section 3,  we first obtain a set of  open-loop control laws in terms of  MF FBSDEs, and give its feedback representation by virtue of Riccati equations.
	In Section 4, we design the feedback strategies of MF Stakelberg games and provide the corresponding costs of all players. In Section 5, we give a numerical example to demonstrate the performance of different solutions. 
	Section 6 concludes the paper.
	
	\emph{Notation}: Throughout this paper, let $(\Omega,\mathcal{F},\{\mathcal{F}_t\}_{0\leq t\leq T}, \mathbb{P})$ be a complete filtered probability space augmented by all $\mathbb{P}$-null sets in $\mathcal{F}$. 
	$|\cdot|$ is the standard Euclidean norm and $\langle\cdot,\cdot\rangle$ is the standard Euclidean inner product.
	For a vector $z$ and a symmetric matrix $Q$, $|z|_Q^2= z^TQz$; $Q>0$ ($Q\geq0$) means that the matrix $Q$ is positive definite (positive semi-definite).
	$Q^{\dag}$ is the Moore-Penrose pseudoinverse\footnote{$Q^{\dag}$ 
		is a unique matrix satisfying
		$QQ^{\dag}Q=Q^{\dag}, Q^{\dag}QQ^{\dag}=Q, (Q^{\dag}Q)^T=Q^{\dag}Q$, and $(QQ^{\dag})^T=QQ^{\dag}$. See \cite{P55} for more properties of pseudoinverse.} of the matrix $Q$, $\mathcal{R}(Q)$ denotes the range of a matrix (or an operator) $Q$.
	Let 
	$C(0,T;\mathbb{R}^{m\times n})$ be the set of $\mathbb{R}^{m\times n}$-valued continuous function and $L_{\mathcal{F}}^2(0,T;\mathbb{R}^m)$  be the set of all $\{\mathcal{F}\}_{t\geq 0}$-adapted $\mathbb{R}^m$-valued processes $x(\cdot)$ such that $\|x(t)\|_{L^2}^2=:\mathbb{E}\int_{0}^{T}\|x(t)\|^2dt<\infty$. 
	
	\section{Problem Formulation}
	\def\theequation{2.\arabic{equation}}
	\setcounter{equation}{0}
	Consider a large-population system with a leader and $N$ followers. The state processes of a leader and  $N$ followers satisfy the following stochastic differential equations: 
	\begin{equation}\label{eq1}
		\left\{
		\begin{aligned}
			dx_0(t)\!=&[A_0x_0(t)+G_0x^{(N)}(t)+B_0u_0(t)]dt+[C_0x_0(t)+\bar{G}_0x^{(N)}(t)+D_0u_0(t)]dW_0(t),\\
			dx_i(t)=&[Ax_i(t)+Bu_i(t)+Gx^{(N)}(t)+Fx_0(t)]dt+[Cx_i(t)+Du_i(t)+\bar{G}x^{(N)}(t)+\bar{F}x_0(t)]dW_i(t),\\
			x_0(0)\!=&\xi_0, \quad x_i(0)=\xi_{i}, \quad i=1,2,\cdots,N,
		\end{aligned}
		\right.
	\end{equation}\\
	where $x_0\in \mathbb{R}^{n_0},u_0\in \mathbb{R}^{m_0}$ are the state and input of the leader, and $x_i\in \mathbb{R}^n,u_i\in \mathbb{R}^m$ are the state and input of the $i$th follower, $i=1,\cdots,N$, respectively. $x^{(N)}(t)\triangleq\frac{1}{N}\sum_{i=1}^{N}x_i(t)$ is the state average of all the followers. $\{W_0(\cdot), W_1(\cdot),\cdots,W_N(\cdot)\}$ are a sequence of independent $d$-dimensional standard Brownian motions defined on the space $(\Omega,\mathcal{F},\{\mathcal{F}_t\}_{0\leq t\leq T}, \mathbb{P})$. 
	Let $\mathcal{F}_t=\sigma(\xi_0, \xi_i, W_0(s), W_i(s), 0\leq s\leq t, i=1,\cdots,N)$). Denote $\mathcal{F}_t^0=\sigma(\xi_0, W_0(s), 0\leq s\leq t)$ and $\mathcal{F}_t^i=\sigma(\xi_0, \xi_i,W_0(s), W_i(s), 0\leq s\leq t)$ for $i=1,\cdots, N$. 
	The admissible control set for the leader is defined as follows:
	$
	\mathcal{U}_{0}=\big\{u_0|u_0(t)\in L_{\mathcal{F}_t^0}^2(0,T;\mathbb{R}^m)\big\}.$
	The admissible 
	decentralized control set for all the followers is defined by
	\begin{equation}\nonumber
		\begin{split}
			\mathcal{U}_{d}=&\Big\{(u_1,\cdots,u_N)|u_i(t) \in L_{\mathcal{F}_t^i}^2(0,T;\mathbb{R}^m), i=1,\cdots,N\Big\}.
		\end{split}
	\end{equation}
	Also, the centralized control set for followers is given by
	\begin{equation*}
		\mathcal{U}_c = \Big\{(u_1,\cdots,u_N)|u_i(t)\in L_{\mathcal{F}_t}^2(0,T;\mathbb{R}^m), i=1,\cdots,N\Big\}.
	\end{equation*}
	
	For the leader, the cost functional is defined by
	\begin{equation}\label{eq2a}
		\begin{aligned}
			{J}_0(u_0,u)=&\mathbb{E}\int_{0}^{T}\big[|x_0(t)- \Gamma_0x^{(N)}(t)|_{Q_0}^{2}+|u_0(t)|_{R_0}^{2}\big]dt+\mathbb{E}\big[|x_0(T)-\hat{\Gamma}_0x^{(N)}(T)|^2_{H_0}\big],
		\end{aligned}
	\end{equation}
	where $Q_0$, $R_0$ and $H_0$ are symmetric  matrices with proper dimensions, and $u=(u_1,\cdots,u_N)$. For the $i$th follower, the cost functional is defined by
	\begin{equation}\label{eq2b}
		\begin{aligned}
			J_i(u_0,u)\!=&\mathbb{E}\int_{0}^{T}\!\!\big[|x_i(t)- \Gamma x^{(N)}(t)-\Gamma_1x_0(t)|_{Q}^{2}+|u_i(t)|_{R}^{2}\big]dt+\mathbb{E}\big[|x_i(T)-\hat{\Gamma}x^{(N)}(T)-\hat{\Gamma}_1x_0(T)|^2_{H}\big],
		\end{aligned}
	\end{equation}
	where $Q$, $R$ and $H$ are  symmetric  matrices with proper dimensions. All the followers cooperate to minimize their social cost functional, denoted by
	\begin{equation}\label{eq4}
		{J}_{\rm soc}^{(N)}(u_0,u)=\frac{1}{N}\sum_{i=1}^{N}{J}_i(u_0,u).
	\end{equation}

	Now we make the following assumption.
	
	(\textbf{A1}) 
	$\{x_i(0)\}$ and $W_i(t),i=1,2,\cdots,N$ are independent of each other. $\mathbb{E}x_0(0)=\bar{\xi}_0$ and
	$\mathbb{E}x_i(0)=\bar{\xi}$, $i=1,\cdots,N$. There exists a constant $c_0$ such that $\sup_{i=1,\cdots,N}\mathbb{E}|x_i(0)|^2\leq c_0$, where $c_0$ is independent of $N$. 
	

	We next discuss the decision hierarchy of the Stackelberg game. The leader 
	holds
	a dominant position in the sense that it first
	announces its strategy $u_0$, and enforces on followers.
	The $N$ followers then respond by cooperatively optimizing their social cost (\ref{eq4})
	under the leader's strategy. In
	this process, the leader takes into account of the rational reactions
	of followers. 

	Due to accessible information restriction and high computational
	complexity, one generally is not able to attain centralized Stackelberg equilibria, but only achieve  asymptotic Stackelberg equilibria under decentralized information patterns.

	We now introduce the definition of the open-loop ($\epsilon_1,\epsilon_2$)-Stackelberg equilibrium.
	From now on, the notation of time $t$ may be suppressed if necessary.
	
	\begin{definition}
		A set of control laws $(u_0^*,
		u_1^*,\cdots,u_N^*)$
		is an open-loop ($\epsilon_1,\epsilon_2$)-Stackelberg equilibrium 
		if the following hold:
		
		(i) When the leader announces a strategy $u_0^*(\cdot)\in \mathcal{U}_0$ over $[0,T]$, $u^*=(u_1^*,\cdots,u^*_N)$ attains an $\epsilon_1$-optimal response, 
		i.e.,
		$J_{\rm soc}^{(N)}(u_0^*,u^*)\leq J_{\rm soc}^{(N)}(u_0^*,u)+\epsilon_1, \ \hbox{for any}\ u\in \mathcal{U}_c, $

		(ii) For any $u_0\in \mathcal{U}_0$, $J_0(u^*_0,u^*(u_0^*))\leq J_0(u_0,u^*(u_0))+\epsilon_2$, where $u^*(\cdot)$ is $\epsilon_1$-optimal response to the leader's strategy.
		
	\end{definition}
	
{Inspired by \cite{BO99, WZ21, YH21}, we consider feedback strategies with the following form:}
	\begin{equation} \label{eq3a}
		\left\{
		\begin{aligned}
			u_0=&P_0x_0+\bar{P}\bar{x},\cr
			u_i=&\hat{K}x_i+\bar{K}\bar{x}+K_0x_0, \ i=1,\cdots,N
		\end{aligned}
		\right.
	\end{equation}
	where $P_0,\bar{P},\hat{K},\bar{K},{K}_0\in L_2(0,T;\mathbb{R}^{n\times n})$; $x_0,x_i$ and $\bar{x}$ satisfy
	\begin{equation}\label{eq3b}
		\left\{
		\begin{aligned}
			&dx_0=[A_0x_0+B_0(P_0x_0+\bar{P}\bar{x})]dt+[C_0x_0+D_0(P_0x_0+\bar{P}\bar{x})]dW_0,\\
			&dx_i=[Ax_i+B(\hat{K}x_i+\bar{K}\bar{x}+K_0x_0)+Gx^{(N)}+Fx_0]dt\cr
			&\hspace{2.5em} +[Cx_i+D(\hat{K}x_i+\bar{K}\bar{x}+K_0x_0)+\bar{G}x^{(N)}+\bar{F}x_0]dW_i,\cr
			&d\bar{x}=\big\{[A+G+B(\hat{K}+\bar{K})]\bar{x}+(F+BK_0)x_0\big\}dt,\cr
			&x_0(0)=\xi_0, \ x_i(0)=\xi_{i}, \ i=1,2,\cdots,N, \ \bar{x}(0)=\bar{\xi}.
		\end{aligned}
		\right.
	\end{equation}
	In the above, $\bar{x}=\mathbb{E}[x_i|\mathcal{F}^0_t]$ is an approximation of $x^{(N)}$ for sufficiently large $N$.

	We now introduce the definition of the feedback ($\epsilon_1,\epsilon_2$)-Stackelberg equilibrium.
	
	\begin{definition}
		A set of strategies $(\hat{u}_0,
		\hat{u}_1,\cdots,\hat{u}_N)$
		is a feedback ($\epsilon_1,\epsilon_2$)-Stackelberg equilibrium 
		if the following hold:
		
		(i)  When the leader announces a strategy $\hat{u}_0=P_0x_0+\bar{P}\bar{x}$ at time $t$, $\hat{u}=(\hat{u}_1,\cdots,\hat{u}_N)$ attains an $\epsilon_1$-optimal feedback response, i.e.,
		$J_{\rm soc}^{(N)}(\hat{u}_0,\hat{u})\leq J_{\rm soc}^{(N)}(\hat{u}_0,u)+\epsilon_1, \ \hbox{for any}\ u\in \mathcal{U}_c, $
		where both $\hat{u}_i$ and $u_i$ have the form $\hat{K}x_i+\bar{K}\bar{x}+K_0x_0$, $i=1,\ldots N$;
		
		(ii) For any $u_0\in \mathcal{U}_0$, $J_0 (\hat{u}_0,\hat{u}(\hat{u}_0) )\leq J_0(u_0,\hat{u}({u}_0))+\epsilon_2$, where $u_0$ has the form $P_0x_0+\bar{P}\bar{x}$, and $\hat{u}(\cdot)$ is $\epsilon_1$-optimal feedback response to the leader's strategy.
		
	\end{definition}
	
	
	%
	{ In this paper, we study open-loop and feedback  solutions to Problem (\ref{eq1})-(\ref{eq4}), respectively.

\textbf{(PO)} Seek an open-loop ($\epsilon_1,\epsilon_2$)-Stackelberg equilibrium over  decentralized control sets $ \mathcal{U}_0$, $ \mathcal{U}_d$;

\textbf{(PF)} Seek a feedback ($\epsilon_1,\epsilon_2$)-Stackelberg equilibrium in  the form of (\ref{eq3a}). }
	
	\section{Open-loop Solutions to Leader-Follower MF Games}
	\def\theequation{3.\arabic{equation}}
	\setcounter{equation}{0}
	\subsection{The MF Social Control Problem for $N$ Followers}
	
	Denote $$\begin{aligned}
		&Q_{\Gamma}\stackrel{\Delta}{=}Q\Gamma+\Gamma^TQ-\Gamma^TQ\Gamma, \ H_{\hat{\Gamma}}\stackrel{\Delta}{=}H\hat{\Gamma}+\hat{\Gamma}^TH-\hat{\Gamma}^TH\hat{\Gamma},\\
		& Q_{\Gamma_1}\stackrel{\Delta}{=}(I-\Gamma)^TQ\Gamma_1,\ H_{\hat{\Gamma}_1}\stackrel{\Delta}{=}(I-\hat{\Gamma})^TH\hat{\Gamma}_1.
	\end{aligned}$$
	
	Suppose $u_0$ 
	is fixed.  
	We now consider the following social control problem for $N$ followers.
	
	\textbf{(P1)}: minimize $ {J}_{\rm soc}$ over $u\in{\mathcal U}_{c}$, where
	{\begin{align*}
			{J}_{\rm soc}(u)=&\sum_{i=1}^N
			\mathbb{E}\int_0^{T}
			\Big[\big|x_i
			-\Gamma {x}^{(N)}-\Gamma_1 x_0\big|^2_{Q}
			+|u_i|^2_{R}\Big]dt+\frac{1}{N}\sum_{i=1}^N\mathbb{E}\big[|x_i(T)-\hat{\Gamma}{x}^{(N)}(T)-\hat{\Gamma}_1x_0(T)|^2_{H}\big].
	\end{align*}}
	
	By examining the social cost variation, we obtain the optimal control laws for $N$ followers.
	\begin{theorem}\label{thm3.1}
		Problem (P1) admits an optimal control if and only if $ {J}_{\rm soc}$ is convex in $u$ and
		the following system of FBSDEs admits a set of adapted solutions $\{x_i,p_i,q_i^j,i,j=1,\cdots,N\}$:
		\begin{equation}\label{eq3}
			\left\{
			\begin{aligned}
dx_0=& (A_0x_0+B_0\check{u}_0+G_0x^{(N)})dt+(C_0x_0+D_0\check{u}_0+\bar{G}_0x^{(N)})dW_0,\ x_0(0)=\xi_0,\cr
dx_i= &(Ax_i+B\check{u}_i+Gx^{(N)}+Fx_0)dt+(Cx_i+D\check{u}_i+\bar{G}x^{(N)}+\bar{F}x_0)dW_i,\ x_i(0)=\xi_i,\cr
	dp_0=& -\big[A_0^Tp_0+F^Tp^{(N)}+C_0^Tq_0^0+\bar{F}^Tq^{(N)}-Q_{\Gamma_1}^T x^{(N)}+\Gamma_1^TQ\Gamma_1 x_0)\big]dt+\sum_{j=0}^Nq_0^jdW_j,\cr
dp_i=&\!-\big(A^Tp_i+G^Tp^{(N)}+G_0^Tp_0+C^Tq_i^i+\bar{G}^Tq^{(N)}+\bar{G}_0^Tq^0_0+Qx_i-Q_{\Gamma}x^{(N)}
				-Q_{\Gamma_1} x_0\big)dt\!+\!\sum_{j=0}^Nq_i^jdW_j,\cr
p_0(T)	&\!=-H_{\bar{\Gamma}_1} ^T{x}^{(N)}(T)+\bar{\Gamma}_1^TH\bar{\Gamma}_1{x}_0(T),\ p_i(T)=Hx_i(T)-H_{\hat{\Gamma}}x^{(N)}(T)-H_{\hat{\Gamma}_1} x_0(T),
			\end{aligned}\right.
		\end{equation}
		where $p^{(N)}=\frac{1}{N}\sum_{j=1}^Np_j$, $q^{(N)}=\frac{1}{N}\sum_{j=1}^Nq_j^j$, and
		the optimal control laws of followers $\check{u}_i$ satisfy 
		\begin{equation}  \label{eq3-b}
 R\check{u}_i+B^Tp_i+D^Tq_i^i=0,\ i=1,\cdots,N.
		\end{equation}
			\end{theorem}
	
	\emph{Proof.} See Appendix A. 
\hfill{$\Box$}

	 {The above theorem gives an equivalence between the solvability of Problem (P1) and that of an FBSDE under the convexity assumption. We refer to the backward equation in (3.2) as the adjoint equation of (1.1). Condition (\ref{eq3-b}) can be regarded as the stationarity condition in Pontryagin's maximum principle.} Indeed, if $J_{\rm soc}^{(N)}$ is uniformly convex in $u$, then Problem (P1) admits an optimal control necessarily \cite{YZ99}. For further existence analysis, we assume
	
	(\textbf{A2}) $J^{(N)}_{\rm soc}$ is uniformly convex in $u$.

\begin{remark}
{ The uniform convexity of $J_{\rm soc}^{(N)}$ in Problem (P1)  can be verified by virtue of the solvability of Riccati equations (See e.g., \cite{SLY16}, \cite{WZ21}).
  Particularly, if $Q\geq0$ and $R>0$, then A2) holds.}
\end{remark}
	
	Denote $\mathbb{E}_{\mathcal{F}^0}[\cdot]\stackrel{\Delta}{=}\mathbb{E}[\cdot|\mathcal{F}^0_t]$. 
	Letting $N\to\infty$, by the MF methodology \cite{hcm2012}, \cite{CD18}, \cite{WZ21}, we can approximate $\check{x}_i$, $\check{p}_i$  in (\ref{eq3}) by $\bar{x}_i$, $\bar{p}_i$, $i=1,\cdots,N$, which satisfy
	\begin{equation}\label{eq5a}
		\left\{
		\begin{aligned}
d\bar{x}_0=& (A_0\bar{x}_0+B_0\check{u}_0+G_0\mathbb{E}_{\mathcal{F}^0}[\bar{x}_i])dt+(C_0\bar{x}_0+D_0\check{u}_0+\bar{G}_0\mathbb{E}_{\mathcal{F}^0}[\bar{x}_i])dW_0,\ \bar{x}_0(0)=\xi_{0},\cr
d\bar{x}_i=& (A\bar{x}_i+B{u}_i^*+G\mathbb{E}_{\mathcal{F}^0}[\bar{x}_i]+F\bar{x}_0)dt+(C\bar{x}_i+D{u}_i^*+\bar{G}\mathbb{E}_{\mathcal{F}^0}[\bar{x}_i]+\bar{F}\bar{x}_0)dW_i,\ \bar{x}_i(0)=\xi_{i},\cr
d\bar{p}_0=& -\big[A_0^T\bar{p}_0+F^T\mathbb{E}_{\mathcal{F}^0}[\bar{p}_i]+C_0^T\bar{q}_0^0+\bar{F}^T\mathbb{E}_{\mathcal{F}^0}[\bar{q}_i^i]-Q_{\Gamma_1}^T \mathbb{E}_{\mathcal{F}^0}[\bar{x}_i]+\Gamma_1^TQ\Gamma_1 \bar{x}_0\big]dt+\bar{q}_0^0dW_0,\cr
&p_0(T)=-H_{\bar{\Gamma}_1} ^T\mathbb{E}_{\mathcal{F}^0}[\bar{x}_i(T)](T)+\bar{\Gamma}_1^TH\bar{\Gamma}_1\bar{x}_0(T),\cr
d\bar{p}_i=&-\big(A^T\bar{p}_i+G^T\mathbb{E}_{\mathcal{F}^0}[\bar{p}_i]+G_0^T\bar{p}_0+C^T\bar{q}_i^i+\bar{G}^T\mathbb{E}_{\mathcal{F}^0}[\bar{q}_i^i]+\bar{G}_0^T\bar{q}_0^0+Q\bar{x}_i
			-Q_{\Gamma}\mathbb{E}_{\mathcal{F}^0}[\bar{x}_i]\cr
			&-Q_{{\Gamma}_1} \bar{x}_0\big)dt+\bar{q}_i^0dW_0+\bar{q}_i^idW_i,\ \bar{p}_i(T)=H\bar{x}_i(T)-H_{\hat{\Gamma}}\mathbb{E}_{\mathcal{F}^0}[\bar{x}_i(T)]-H_{\hat{\Gamma}_1} \bar{x}_0(T),
		\end{aligned}\right.
	\end{equation}
	with the decentralized control ${u}_i^*$ {satisfying the stationarity condition}
	\begin{equation}\label{eq6-b}
		R{u}_i^*+B^T\bar{p}_i+D^T\bar{q}_i^i=0, \ i=1,\cdots,N.
	\end{equation}
	
	We now use the idea inspired by {\cite{my1999}, \cite{Y13}, \cite{WZZ22}} to decouple the FBSDE
	(\ref{eq5a}).
	Let $\bar{p}_0=K\bar{x}_0+\bar{K}\mathbb{E}_{\mathcal{F}^0}[\bar{x}_i]+\varphi_0$ and $\bar{p}_i=P\bar{x}_i+\bar{P}\mathbb{E}_{\mathcal{F}^0}[\bar{x}_i]+P_0\bar{x}_0+\varphi,\ i=1,\cdots,N,$ where $\varphi_0$ satisfies $d\varphi_0=\hat{\varphi}_0dt+\zeta_0dW_0$ and $\varphi$ satisfies $d\varphi=\hat{\varphi}dt+\zeta dW_0$.
 Then, we have
\begin{equation}\label{eq10a}
		\begin{aligned}
  	d\bar{p}_0=&\dot{K}\bar{x}_0dt+K\big[ (A_0\bar{x}_0+B_0\check{u}_0+G_0\mathbb{E}_{\mathcal{F}^0}[\bar{x}_i])dt+(C_0\bar{x}_0+D_0\check{u}_0+\bar{G}_0\mathbb{E}_{\mathcal{F}^0}[\bar{x}_i])dW_0  \big]\cr
  &+\dot{\bar{K}}\mathbb{E}_{\mathcal{F}^0}[\bar{x}_i]dt+\bar{K}\big[(A+G)\mathbb{E}_{\mathcal{F}^0}[\bar{x}_i]+B\mathbb{E}_{\mathcal{F}^0}[\bar{u}_i]+F\bar{x}_0\big]+d{\varphi}_0\cr
  =& -\big[A_0^T(K\bar{x}_0+\bar{K}\mathbb{E}_{\mathcal{F}^0}[\bar{x}_i]+\varphi_0)+F^T\big((P+\bar{P})\mathbb{E}_{\mathcal{F}^0}[\bar{x}_i]+P_0\bar{x}+\varphi\big)+C_0^T\bar{q}_0^0\cr
  &+\bar{F}^T\mathbb{E}_{\mathcal{F}^0}[\bar{q}_i^i]-Q_{\Gamma_1}^T \mathbb{E}_{\mathcal{F}^0}[\bar{x}_i]+\Gamma_1^TQ\Gamma_1 \bar{x}_0\big]dt+\bar{q}_0^0dW_0,
\end{aligned}
	\end{equation}
which implies 	\begin{equation}\label{eq10a1}
  \bar{q}_0^0=K(C_0\bar{x}_0+D_0\check{u}_0+\bar{G}_0\mathbb{E}_{\mathcal{F}^0}[\bar{x}_i])+\zeta_0,
\end{equation}
and
	\begin{equation}\label{eq10}
		\begin{aligned}
			d\bar{p}_i=\ &\dot{P}\bar{x}_idt+P\Big[\big(A\bar{x}_i+B{u}_i^*+G\mathbb{E}_{\mathcal{F}^0}[\bar{x}_i]+F\bar{x}_0\big)dt+(C\bar{x}_i+D{u}_i^*+\bar{G}\mathbb{E}_{\mathcal{F}^0}[\bar{x}_i]+\bar{F}x_0) dW_i\Big]\cr
			&+\dot{\bar{P}}\mathbb{E}_{\mathcal{F}^0}[\bar{x}_i]dt+\bar{P}\big[(A+G)\mathbb{E}_{\mathcal{F}^0}[\bar{x}_i]+B\mathbb{E}_{\mathcal{F}^0}[\bar{u}_i]+F\bar{x}_0\big]dt+d{\varphi}\cr
	&+\dot{P}_0\bar{x}_0dt+P_0\big[(A_0\bar{x}_0+B_0\check{u}_0+G_0\mathbb{E}_{\mathcal{F}^0}[\bar{x}_i])dt+(C_0\bar{x}_0+D_0\check{u}_0+\bar{G}_0\mathbb{E}_{\mathcal{F}^0}[\bar{x}_i])dW_0\big]\cr		
=\ &-\Big[A^T(P\bar{x}_i+\bar{P}\mathbb{E}_{\mathcal{F}^0}[\bar{x}_i]+P_0\bar{x}_0+\varphi)+C^T [\bar{q}_i^i]+G^T((P+\bar{P})\mathbb{E}_{\mathcal{F}^0}[\bar{x}_i]+P_0\bar{x}_0+\varphi)\cr
			&+G_0^T(K\bar{x}_0+\bar{K}\mathbb{E}_{\mathcal{F}^0}[\bar{x}_i]+\varphi_0)+\bar{G}^T\mathbb{E}_{\mathcal{F}^0}[\bar{q}_i^i]+\bar{G}_0^T\bar{q}_0^0
\cr			&+Q\bar{x}_i-Q_{\Gamma} \mathbb{E}_{\mathcal{F}^0}[\bar{x}_i]-Q_{\Gamma_1} \bar{x}_0\Big]dt+\bar{q}_i^idW_i+\bar{q}_i^0dW_0,
		\end{aligned}
	\end{equation}
	which implies
	\begin{equation}\label{eq14a}
		\bar{q}_i^i=P(C\bar{x}_i+D{u}_i^*+\bar{G}\mathbb{E}_{\mathcal{F}^0}[\bar{x}_i]+\bar{F}\bar{x}_0), \ i=1,\cdots,N.
	\end{equation}
	This together with (\ref{eq6-b}) leads to
	\begin{align*}
		&R{u}_i^*+B^T(P\bar{x}_i+\bar{P}\mathbb{E}_{\mathcal{F}^0}[\bar{x}_i]+P_0\bar{x}_0+\varphi)+D^TP(C\bar{x}_i+D\bar{u}_i+\bar{G}\mathbb{E}_{\mathcal{F}^0}[\bar{x}_i]+\bar{F}\bar{x}_0)=0.
	\end{align*}
	Let $\Upsilon\stackrel{\Delta}{=}R+D^TPD$. If $\mathcal{R}(B^T)\cup  \mathcal{R}(D^TP)\subseteq  \mathcal{R}(\Upsilon)$, then we have
	\begin{align}\label{eq10e}
		{u}^*_i=&-\Upsilon^{\dag}\big[\big(B^TP+D^TPC\big)\bar{x}_i+\big(B^T\bar{P}+D^TP\bar{G}\big)\mathbb{E}_{\mathcal{F}^0}[\bar{x}_i]+B^T\varphi+(B^TP_0+D^TP\bar{F})\bar{x}_0\big].
	\end{align}
From (\ref{eq10a}), (\ref{eq10a1}) and (\ref{eq14a}), the following equations should hold:
	\begin{equation}\label{eq9b}
  \left\{\begin{aligned}
  &\dot{K}+K^{T}A_0+A^T_0K+C_0^TKC_0+\bar{F}^TP\bar{F}-\Psi_1^T\Upsilon^{\dag}(B^TP_0+D^TP\bar{F})\cr
  &+\bar{K}F+F^TP_0+\Gamma_1^TQ\Gamma_1=0, \ K(T)=\hat{\Gamma}_1^TH\hat{\Gamma}_1,\cr
   & \dot{\bar{K}}+\bar{K}(A+G)+A_0^T\bar{K}-\Psi_1^T\Upsilon^{\dag}\big[B^T(P+\bar{P})+D^TP(C+\bar{G})\big]+F^T(P+\bar{P})\cr
   &+C_0^TK\bar{G}_0
    +\bar{F}^TP(C+\bar{G})+KG_0-Q^T_{\Gamma_1}=0,\ \bar{K}(T)=H^T_{\hat{\Gamma}_1},\cr
    & d\varphi_0+A_0^T\varphi_0+(F-B\Upsilon^{\dag}\Psi_1 )^T\varphi+C_0^T\zeta_0+(C^T_0KD_0+KB_0)\check{u}_0-\zeta_0dW_0=0,  \varphi_0(T)=0,
\end{aligned}\right.
\end{equation}
where $\Psi_1\stackrel{\Delta}{=}B^T\bar{K}^T+D^TP\bar{F}$.
	From (\ref{eq10a1})-(\ref{eq10e}),   the following should hold:
		\begin{equation} \label{eq10b}
\left\{\begin{aligned}
		&\dot{P}+A^TP+PA+C^TPC+Q-\Psi_2^T\Upsilon^{\dag}\Psi_2=0, P(T)=H,\\
		&\dot{\bar{P}}+(A+G)^T\bar{P}+\bar{P}(A+G)+G^TP+PG+P_0G_0+G_0^T\bar{K}+C^TP\bar{G}+\bar{G}^TP(C+\bar{G})\\
		&\ \ \, +\bar{G}_0^TK\bar{G}_0-\Psi_2^T\Upsilon^{\dag}\Psi_3-\Psi_3^T\Upsilon^{\dag}\Psi_2-\Psi_3^T\Upsilon^{\dag}\Psi_3-Q_{\Gamma}=0,\ \bar{P}(T)=-H_{\hat{\Gamma}},\\
		&\dot{P}_0+P_0A_0+(A+G)^TP_0+(C+\bar{G})^TP\bar{F}-(\Psi_2+\Psi_3)^T\Upsilon^{\dag}(B^TP_0+D^TP\bar{F})\cr
&+(P+\bar{P})F+\bar{G}_0KC_0+G_0^T{K}-Q_{\Gamma_1}=0,\ P_0(T)=H_{\hat{\Gamma}_1},\cr
&d{\varphi}+\big\{\big[A+G-B\Upsilon^{\dag}\big(\Psi_2+\Psi_3)\big)\big]^T\varphi+G_0^T\varphi_0+\bar{G}_0^T\zeta_0+P_0B_0\check{u}_0\big\}dt
-\zeta dW_0=0,\ \varphi(T)=0,
	\end{aligned}\right.
\end{equation}
where $\Psi_2=B^TP+D^TPC$ and $\Psi_3=B^T\bar{P}+D^TP\bar{G}$.
It can be verified that $\bar{K}^T=P_0$, and matrices $K,P,\bar{P}$ are symmetric.
We assume
	
	\textbf{(A3)} Equations (\ref{eq9b})-(\ref{eq10b}) admit a set of solution ($K,\bar{K},P,\bar{P},P_0,\varphi_0,\varphi$) such that $\Upsilon\geq 0$, and
	\begin{equation*}
		\mathcal{R}(B^T)\cup  \mathcal{R}(D^TP)\subseteq  \mathcal{R}(\Upsilon).
	\end{equation*}
	

	From the above discussion, we have the following result.
	\begin{proposition}
		Under (A3),  the decentralized control in (\ref{eq6-b}) has a feedback representation (\ref{eq10e}).
	\end{proposition}

	Applying (\ref{eq10e}) into (\ref{eq5a}), we obtain that $\bar{x}=\mathbb{E}_{\mathcal{F}^0}[\bar{x}_i]$ satisfies
	\begin{equation}\label{eq7}
		d\bar{x}= \big[\big(A+G-B\Upsilon^{\dag}(\Psi_2+\Psi_3)\big)\bar{x}-B\Upsilon^{\dag}B^T\varphi+(F-B\Upsilon^{\dag}\Psi_1)\bar{x}_0\big]dt.
	\end{equation}
	
	\subsection{Optimization for the Leader}
	
	Denote $\bar{A}\stackrel{\Delta}{=}A-B\Upsilon^{\dag}(B^TP+D^TPC)$, and $\bar{C}\stackrel{\Delta}{=}C-D\Upsilon^{\dag}(B^TP+D^TPC)$.
	After applying the control laws of followers in (\ref{eq10e}), we have the following optimal control problem for the leader.
	
	\textbf{(P2)}: minimize $ {J}_{0}(u_0,u^*(u_0))$ over $u_0\in \mathcal{U}_0$, where
	\begin{equation} \label{eq23-a}
	\left\{\begin{aligned}
		&{J}_0(u_0,u^*(u_0))=\mathbb{E}\int_{0}^{T}\big[|x_0- \Gamma_0{x}_*^{(N)}|_{Q_0}^{2}+|u_0|_{R_0}^{2}\big]dt+\mathbb{E}\big[|x_0(T)-\hat{\Gamma}_0x_*^{(N)}(T)|^2_{H_0}\big],\cr
		&dx_0=(A_0x_0+G_0x^{(N)}+B_0u_0)dt+(C_0x_0+\bar{G}_0x^{(N)}+D_0u_0)dW_0,\ x_0(0)=\xi_0,\cr
		&d{x}_i^*=\big[{A}{x}_i^*+G{x}^{(N)}_*-B\Upsilon^{\dag}\big(\Psi_2\bar{x}_i+\Psi_3\bar{x}+B^T\varphi\big)+Fx_0-B\Upsilon^{\dag}\Psi_1\bar{x}_0\big]dt\\
		&\hspace{3em}+\big[Cx_i^*+\bar{G}{x}^{(N)}_*-D\Upsilon^{\dag}\big(\Psi_2\bar{x}_i+\Psi_3\bar{x}+B\varphi\big)+\bar{F}x_0-D\Upsilon^{\dag}\Psi_1\bar{x}_0\big]dW_i,\
{x}_i^*(0)=\xi_i,\cr
		&d\varphi_0=-\big[A_0^T\varphi_0+(F-B\Upsilon^{\dag}\Psi_1 )^T\varphi+C_0^T\zeta_0+(C^T_0KD_0+KB_0)\check{u}_0\big]+\zeta_0dW_0=0,\ \varphi_0(T)=0,\\
		&d{\varphi}=-\big\{\big[A+G-B\Upsilon^{\dag}\big(\Psi_2+\Psi_3)\big)\big]^T\varphi+G_0^T\varphi_0+\bar{G}_0^T\zeta_0+P_0B_0\check{u}_0\big\}dt
+\zeta dW_0=0,\ \varphi(T)=0,
	\end{aligned}\right.
	\end{equation}
	where $x_i^*$ is the realized state under control $u^*_i, i=1,\cdots,N$, and $x_*^{(N)}=\frac{1}{N}\sum_{i=1}^Nx_i^*$.
	From (\ref{eq23-a}), we have
	\begin{align*}
		d{x}_*^{(N)}=&\big[({A}+G){x}_*^{(N)}-B\Upsilon^{\dag}\big(\Psi_2\bar{x}^{(N)}+\Psi_3\bar{x}+B^T\varphi\big)+Fx_0-B\Upsilon^{\dag}\Psi_1\bar{x}_0\big]dt\cr
		&+\frac{1}{N}\sum_{i=1}^N\big[Cx_i^*+\bar{G}{x}^{(N)}_*-D\Upsilon^{\dag}\big(\Psi_2\bar{x}_i+\Psi_3\bar{x}+B\varphi\big)+Fx_0-D\Upsilon^{\dag}\Psi_1\bar{x}_0\big]dW_i,\ {x}_*^{(N)}(0)=\frac{1}{N}\sum_{i=1}^N\xi_i,
	\end{align*}
	where $\bar{x}^{(N)}=\frac{1}{N}\sum_{i=1}^N\bar{x}_i$.
	Note that  $\{W_i\}$ are independent Wiener processes and $\{x_i(0)\}$ are independent r.v.s. For the large population case, it is plausible to replace $\bar{x}^{(N)}$, ${x}_*^{(N)}$ by $\bar{x}$, which evolves from (\ref{eq7}). 
	Then we have the limiting optimal control problem for the leader.
	
	\textbf{(P2$^\prime$)}: minimize $ \bar{J}_{0}(u_0,u^*(u_0))$ over $u_0\in {\mathcal{U}}_0$, where
	\begin{align} \label{eq8-a}
		\bar{J}_0(u_0,u^*(u_0))=&\mathbb{E}\int_{0}^{T}\big[|\bar{x}_0- \Gamma_0\bar{x}|_{Q_0}^{2}+|u_0|_{R_0}^{2}\big]dt+\mathbb{E}\big[|\bar{x}_0-\hat{\Gamma}_0\bar{x}(T)|^2_{H_0}\big],
	\end{align}
	subject to
	\begin{equation}\label{eq8c} 
		\left\{
		\begin{aligned}
			&d\bar{x}_0=(A_0x_0+G_0\bar{x}+B_0u_0)dt+(C_0\bar{x}_0+\bar{G}_0\bar{x}+D_0u_0)dW_0,\ \bar{x}_0(0)=\xi_0,\\
			& d\bar{x}=\big[\big(A+G-B\Upsilon^{\dag}(\Psi_2+\Psi_3)\big)\bar{x}-B\Upsilon^{\dag}B^T\varphi+(F-B\Upsilon^{\dag}\Psi_1)\bar{x}_0\big]dt,\\
		&d\varphi_0=-\big[A_0^T\varphi_0+(F-B\Upsilon^{\dag}\Psi_1 )^T\varphi+C_0^T\zeta_0+(C^T_0KD_0+KB_0){u}_0\big]+\zeta_0dW_0=0,\ \varphi_0(T)=0,\\
		&d{\varphi}=-\big\{\big[A+G-B\Upsilon^{\dag}\big(\Psi_2+\Psi_3)\big)\big]^T\varphi+G_0^T\varphi_0+\bar{G}_0^T\zeta_0+P_0B_0{u}_0\big\}dt
+\zeta dW_0=0,\ \varphi(T)=0.
		\end{aligned}\right.
	\end{equation}
	
	We first provide the condition under which Problem (P2$^\prime$) is convex.
	The proof is similar to \cite{jm2017}, \cite{WZ21}, and so omitted here.
	
	\begin{lemma}\label{lem4.1}
		$\bar{J}_0(u_0,{u}^*(u_0))$ is convex in $u_0$ if and only if $\bar{J}_0^0(u_0,{u}^*(u_0))\geq 0$, where
		\begin{align*}
			\bar{J}_0^0(u_0,u^*)=&\mathbb{E}\int_{0}^{T}\big[|\bar{x}_0^0- \Gamma_0\bar{x}^0|_{Q_0}^{2}+|u_0|_{R_0}^{2}\big]dt+\mathbb{E}\big[|\bar{x}_0^0(T)-\hat{\Gamma}_0\bar{x}^0(T)|^2_{H_0}\big],
		\end{align*}
		subject to
		\begin{equation}\label{eq4.5}
			\left\{\begin{aligned}
				d\bar{x}_0^0=& (A_0x_0^0+G_0\bar{x}^0+B_0u_0)dt+(C_0x_0^0+\bar{G}_0\bar{x}^0+D_0u_0)dW_0,\ \bar{x}_0^0(0)=0,\\
				d\bar{x}^0=&\big[\big(A+G-B\Upsilon^{\dag}(\Psi_2+\Psi_3)\big)\bar{x}^0-B\Upsilon^{\dag}B^T\varphi^0+(F-B\Upsilon^{\dag}\Psi_1)\bar{x}_0^0\big]dt,\ \bar{x}^0(0)=0,\\
d\varphi_0^0=&-\big[A_0^T\varphi_0^0+(F-B\Upsilon^{\dag}\Psi_1 )^T\varphi^0+C_0^T\zeta_0^0+(C^T_0KD_0+KB_0){u}_0\big]+\zeta_0^0dW_0=0,\ \varphi_0^0(T)=0,\\
d{\varphi}^0=&-\big\{\big[A+G-B\Upsilon^{\dag}\big(\Psi_2+\Psi_3)\big]^T\varphi^0+G_0^T\varphi_0^0+\bar{G}_0^T\zeta_0^0+P_0B_0{u}_0\big\}dt
+\zeta^0 dW_0=0,\ \varphi^0(T)=0.
			\end{aligned}\right.
		\end{equation}
	\end{lemma}
	
	We now give the following maximum principle for  (P2$^{\prime}$).
	
	\begin{theorem}\label{thm4.0}
		Assume (A1)-(A3) hold.
		Problem (P2$^{\prime}$) admits an 
		optimal control $u_0^*$ if and only if  $\bar{J}_0(u_0,{u}^*(u_0))$ is convex in $u_0$, and 
		the following FBSDE
	\begin{equation}\label{eq16}
			\left\{
			\begin{aligned}
				dy_0=&-\big\{ A^T_0 y_0+C_0^T\beta_0+(F-B\Upsilon^{\dag}\Psi_1)^T\bar{y}+Q_0(\bar{x}^*_0-\Gamma_0\bar{x}^*)\big\}dt+\beta_0dW_0,\\
				d\bar{y}=&-[\big(A+G-B\Upsilon^{\dag}(\Psi_2+\Psi_3)\big)^T\bar{y}+G_0^Ty_0+\bar{G}_0^T\beta_0-\Gamma_0^TQ_0(x_0^*-\Gamma_0\bar{x}^*)]+\bar{\beta}dW_0,\\
	d\psi_0=&\ (A_0\psi_0+G_0\psi)dt+(C_0\psi_0+\bar{G}_0\psi)dW_0,\ \\
				d\psi=&\big[\big(A+G-B\Upsilon^{\dag}(\Psi_2+\Psi_3)\big)^T\psi+(F-B\Upsilon^{\dag}\Psi_1 )\psi_0+B\Upsilon^{\dag}B^T \bar{y}\big]dt,\\
y_0(T)&=H_0(\bar{x}_0(T)-\hat{\Gamma}_0 \bar{x}^*(T)),\ \bar{y}(T)=-\hat{\Gamma}_0^TH_0(\bar{x}_0^*(T)-\hat{\Gamma}_0 \bar{x}^*(T)),\ \psi_0(0)=0,\ \psi(0)=0
			\end{aligned}
			\right.
		\end{equation}
		has a solution such that $u_0^*$ satisfies $R_0u_0^*+B_0^Ty_0+D_0^T\beta_0-(C^T_0KD_0+KB_0)^T\psi_0-(P_0B_0)^T\psi=0$. 

	\end{theorem}
	
	\emph{Proof.}
	Suppose $\{u_0^*\}$ is a candidate of  the optimal control of Problem (P2$^\prime$).
	Let $\bar{x}_0^*$ and $\bar{x}^*$ be the  leader's state and followers' average effect under the control $\{u_0^*\}$.
	Note that
	\begin{align}\label{eq26b}
		\bar{J}_{0}({u}_0^*+\theta u_0,{u}({u}_0^*+\theta u_0)) - \bar{J}_{0}({u}_0^*,{u}^*({u}_0^*))=2\theta I_1+\theta^2I_2,
	\end{align}
	where
	\begin{align}\label{eq29a}
		I_1=& \mathbb{E}\int_0^T\big[\langle Q_0(\bar{x}_0^*-\Gamma_0\bar{x}_*), \bar{x}_0^0-\Gamma_0 \bar{x}^0\rangle+\langle u_0^*, R_0 u_0\rangle\big]dt\\
		&+\mathbb{E}\big[\langle H_0 (\bar{x}_0^*(T)-\hat{\Gamma}_0 \bar{x}^*(T)), \bar{x}_0^0(T)-\hat{\Gamma}_0\bar{x}^0(T)\rangle\big], \label{eq29b} \cr
		I_2=&\mathbb{E}\int_0^T\big[| \bar{x}_0^0-\Gamma_0 \bar{x}^0|^2_{Q_0}+| u_0|^2_{R_0}\big]dt+\mathbb{E}\big[| \bar{x}_0^0(T)-\hat{\Gamma}_0\bar{x}^0(T)|^2_{H_0}\big].
	\end{align}
	Note that for the given $x_0^*$ and $\bar{x}^*$, FBSDE (\ref{eq16}) admits a unique solution. 
	From (\ref{eq4.5}) and (\ref{eq16}), applying It\^{o}'s formula, we obtain
	 \begin{align}\label{eq12}
		&\mathbb{E}[\langle H_0( \bar{x}_0^*-\hat{\Gamma}_0 \bar{x}^*)+\hat{\Gamma}_1^T H(\hat{\Gamma}-I)\psi (T), \bar{x}_0^0(T)\rangle]=\mathbb{E}[\langle y_0(T),\bar{x}_0^0(T)\rangle-\langle y_0(0),\bar{x}_0^0(0)\rangle]\\
		=&\mathbb{E}\int_0^T\Big\{\big\langle -\big[(F-B\Upsilon^{\dag}\Psi_1)^T\bar{y}+Q_0(x_0^*-\Gamma_0\bar{x}^*)\big], \bar{x}_0^0\big\rangle
+\langle G_0^Ty_0+\bar{G}_0^T\beta_0,\bar{x}_0\rangle+\langle B_0^Ty_0+D_0^T\beta_0, u_0\rangle
		\Big\}dt,\label{eq13} \cr
		&-\mathbb{E}[\langle \hat{\Gamma}_0^T H_0(\bar{x}_0^*-\hat{\Gamma}_0 \bar{x}^*), \bar{x}^0(T)\rangle]
=\mathbb{E}[\langle \bar{y}(T), \bar{x}^0(T)\rangle-\langle \bar{y}(0), \bar{x}^0(0)\rangle]\\
		=&\mathbb{E}\int_0^T\big[\langle \Gamma_0^TQ_0(\bar{x}_0-\Gamma_0\bar{x})-G_0^Ty_0-\bar{G}_0^T\beta_0, \bar{x}^0\rangle-\langle B\Upsilon^{\dag}B^T \bar{y}, \varphi^0\rangle+\langle F-B\Upsilon^{\dag}\Psi_1\bar{y}, x_0^0\rangle\big]dt.\nonumber
	\end{align}
	\begin{align}\label{eq13b}
0=&\mathbb{E}[\langle\varphi^0_0(T),\psi_0(T)\rangle -\langle\varphi_0^0(0),\psi_0(0)\rangle]\cr
		=& \mathbb{E}\!\int_0^T\!\!\Big[\langle G_0\psi, \varphi^0_0\rangle-\big\langle (F-B\Upsilon^{\dag}\Psi_1 )^T\varphi^0+C_0^T\zeta_0^0+(C^T_0KD_0+KB_0){u}_0\big], \psi_0\rangle+\langle\zeta_0^0,v_0\rangle\Big] dt\cr
=& \mathbb{E}\!\int_0^T\!\!\Big[\langle G_0\psi, \varphi^0_0\rangle-\big\langle (F-B\Upsilon^{\dag}\Psi_1 )\psi_0,\varphi^0\rangle+\langle v_0-C_0\psi_0,\zeta_0^0\rangle-\langle(C^T_0KD_0+KB_0)^T\psi_0,{u}_0\rangle\Big] dt
	\end{align}
	and
	\begin{align}\label{eq14}
0=&\mathbb{E}[\langle\varphi^0(T),\psi(T)\rangle -\langle\varphi^0(0),\psi(0)\rangle]\\
		=& \mathbb{E}\!\int_0^T\!\!\Big[\langle (F-B\Upsilon^{\dag} \Psi_1)\varphi_0+B\Upsilon^{\dag}B^T \bar{y}, \varphi^0\rangle-\big\langle G_0^T\varphi_0^0+\bar{G}_0^T\zeta_0^0+P_0B_0{u}_0, \psi\rangle\Big] dt\cr
=& \mathbb{E}\!\int_0^T\!\!\Big[\langle (F-B\Upsilon^{\dag} \Psi_1)\varphi_0+B\Upsilon^{\dag}B^T \bar{y}, \varphi^0\rangle
-\langle G_0\psi,\varphi_0^0\rangle-\langle \bar{G}_0\psi, \zeta_0^0\rangle-\langle (P_0B_0)^T\psi,{u}_0\rangle\Big] dt.\nonumber
	\end{align}
		From (\ref{eq29a}) and (\ref{eq12})-(\ref{eq14}), it follows that
		$$\begin{aligned}
			I_1=  &\mathbb{E}\int_0^T\big\langle B_0^Ty_0+D_0^T\beta_0-(C^T_0KD_0+KB_0)^T\psi_0-(P_0B_0)^T\psi+Ru_0^*, u_0\big\rangle dt.
		\end{aligned}$$
Note that $\theta$ is arbitrary.  By (\ref{eq26b}), $u_0^*$ is a minimizer of (P2$^{\prime}$) if and only if $I_1=0$ and $I_2\geq0$.  Indeed, if $I_2\geq 0$ does not hold, 
then there exists some $\check{u}_0\in \mathcal{U}_0$ such that $\bar{J}_{0}^0(\check{u}_0, u^*)< 0$. Then we have $\bar{J}_{0}^0(k\check{u}_0, u^*)=k^{2}\bar{J}_{0}^0(\check{u}_0, u^*) \to -\infty$
($k\to \infty$), which implies the minimization problem should be ill-posed.
		Thus, by Lemma \ref{lem4.1}, $u_0^*$ is an optimal control of (P2$^{\prime}$) if and only if  $Ru_0^*+B_0^Ty_0+D_0^T\beta_0-(C^T_0KD_0+KB_0)^T\psi_0-(P_0B_0)^T\psi=0$ and $\bar{J}_0(u_0,{u}(u_0))$ is convex in $u_0$.
		\hfill{$\Box$}

		%
		Let $X=[x_0^T,\bar{x}^T,\psi_0^T,\psi^T ]^T,Y=[y_0^T,\bar{y}^T\!,\varphi_0^T,\varphi^T]^T,$ $Z=[\beta_0^T,\bar{\beta}^T\!,\zeta_0^T,\zeta^T]^T,
		\mathcal{D}_0=[D_0^T,0,0,0]^T$, $\mathcal{B}_0=[B_0^T,0,(C^T_0KD_0+ KB_0)^T,(P_0B_0)^T]^T$, and
						$$\mathcal{A}=\left[\begin{array}{cccc}
					A_0&G_0&0&0\\
					F-B\Upsilon^{\dag}\Psi_1&{A}+{G}-B\Upsilon^{\dag}(\Psi_2+\Psi_3)&0\\
					0&0&A_0&{G}_0\\
0&0&F-B\Upsilon^{\dag}\Psi_1&{A}+{G}-B\Upsilon^{\dag}(\Psi_2+\Psi_3)
				\end{array}\right],
	$$
$$\mathcal{B}=\left[\begin{array}{cccc}
					0&0&0&0\\
					0&0&0&B\Upsilon^{\dag}B^T\\
0&0&0&0\\
					0&-B\Upsilon^{\dag}B^T&0&0
				\end{array}\right],\quad \mathcal{C}_0=\left[\begin{array}{cccc}
	C_0&	\bar{G}_0&0&0\\
	0&0&0&0\\
0&0& C_0&	\bar{G}_0\\
	0&0&0&0
\end{array}\right],$$
$$\mathcal{H}_0=\left[\begin{array}{cccc}
	H_0&-H_0 \hat{\Gamma}_0&0&0\\
	-\hat{\Gamma}_0^TH_0&\Gamma_0^TH_0\Gamma_0&0&0\\
	0&0&0&0\\
0&0&0&0
\end{array}\right],
\quad \mathcal{Q}\!=\!\!\!\left[\!\begin{array}{cccc}
		-Q_0&Q_0\Gamma_0& 0&0\\
		\\
		\Gamma_0^TQ_0&-\Gamma_0^TQ_0\Gamma_0&0&0\\
		0&0&0&0\\
	0&0&0&0
	\end{array}\!\!\right]\!\!.$$

With above notations, we can rewrite (\ref{eq8c}) and (\ref{eq16})  
as
\begin{equation}
	\label{eq17}
	\left\{
	\begin{aligned}
		dX&=(\mathcal{A}X-\mathcal{B}Y+\mathcal{B}_0u_0^*)dt+(\mathcal{C}_0X+\mathcal{D}_0u_0^*)dW_0,\ X(0)=[\xi_0^T, \bar{\xi}^T,0,0]^T,\\
		dY&=(\mathcal{Q}X-\mathcal{A}^TY-\mathcal{C}_0^TZ)dt+ZdW_0,\ Y(T)=\mathcal{H}_0 X(T),
	\end{aligned}
	\right.
\end{equation}
together with the condition
\begin{equation}\label{4.12}
	R_0u_0^*+\mathcal{B}_0^TY+\mathcal{D}_0^TZ=0.
\end{equation}

We now provide a sufficient condition to guarantee the solvability of (\ref{eq17}).
\begin{proposition}
	Denote $\Upsilon_0{=}R_0+\mathcal{D}_0^T\mathcal{P}\mathcal{D}_0$.
	If  the equation
	\begin{align}\label{eq18}
		&\dot{\mathcal{P}}+\mathcal{P}\mathcal{A}+\mathcal{A}^T\mathcal{P}+\mathcal{C}_0^T\mathcal{P}\mathcal{C}_0
		-\mathcal{Q}-\mathcal{P}\mathcal{B}\mathcal{P}-
		(\mathcal{B}_0^T\mathcal{P}+\mathcal{D}_0^T\mathcal{P}\mathcal{C}_0)^T
		\Upsilon_0^{\dag}(\mathcal{B}_0^T\mathcal{P}+\mathcal{D}_0^T\mathcal{P}\mathcal{C}_0)=0,
	\end{align}
	with $\mathcal{P}(T)=\mathcal{H}_0$ has a solution in $[0,T]$, then FBSDE (\ref{eq17}) is solvable.
\end{proposition}
\emph{Proof.} {Let $\bar{Y}=\mathcal{P}X$ and  $\bar{Z}=\mathcal{P}\big[\mathcal{C}_0-\mathcal{D}_0^T\Upsilon_0^{\dag}(\mathcal{B}_0^T\mathcal{P}+\mathcal{D}_0^T\mathcal{P}\mathcal{C}_0)\big]X$,
where $\mathcal{P}$ is a solution to (\ref{eq18}).
Let $u_0=-\Upsilon_0^{\dag}(\mathcal{B}_0^T\mathcal{P}+\mathcal{D}_0^T\mathcal{P}\mathcal{C}_0)X$.  Denote $\tilde{Y}=Y-\bar{Y}$ and $\tilde{Z}=Z-\bar{Z}$.
Then a direct computation shows
\begin{align*}
  d\tilde{Y}=&
  [(PB-A^T)\tilde{Y}-C_0^T\tilde{Z}]dt+\tilde{Z}dW_0, \quad \tilde{Y}(T)=0.
\end{align*}
It is clear that such a backward SDE admits a unique solution $\tilde{Y}=\tilde{Z}=0$ (\cite{my1999}).
Hence,
 ${Y}=\mathcal{P}X$ and  ${Z}=\mathcal{P}\big[\mathcal{C}_0-\mathcal{D}_0^T\Upsilon_0^{\dag}(\mathcal{B}_0^T\mathcal{P}+\mathcal{D}_0^T\mathcal{P}\mathcal{C}_0)\big]X$.
Then FBSDE (\ref{eq17}) admits an adapted solution.}
\hfill$\Box$

\begin{remark}
	Note that matrices $\mathcal{Q}$  and $\mathcal{H}_0$ are symmetric, and $\mathcal{B}$ are nonsymmetric. We find that  (\ref{eq18}) is a nonsymmetric Riccati equation.
	The existence condition of its solution may be referred in \cite{AF03}, \cite{my1999}.
\end{remark}

For further analysis, assume

\textbf{(A4)} Equation (\ref{eq18}) admits a solution in $C[0,T;\mathbb{R}^{3n}]$.

Under (A4), we construct the following decentralized control laws
\begin{equation}\label{eq27a}
	\left\{
	\begin{aligned}
		{u}_0^*=&-\Upsilon_0^{\dag}(\mathcal{B}_0^T\mathcal{P}+\mathcal{D}_0^T\mathcal{P}\mathcal{C}_0)X,\cr
		{u}_i^*=&-\Upsilon^{\dag}\big[\big(B^TP+D^TPC\big)\bar{x}_i+B^T\varphi+D^TP\bar{F}x_0^*+\big(B^TK+D^TP\bar{G}\big)\mathbb{E}_{\mathcal{F}^0}[\bar{x}_i]\big]
	\end{aligned}
	\right.
\end{equation}
where $X$ and $\bar{x}_i$ are given by (\ref{eq17}), (\ref{eq5a}), and $x_0^*$ is the realized state under the control $u_0^*$.

\begin{theorem}\label{thm4.2}
Assume that (A1)-(A4) hold. Then $({u}_0^*,{u}_1^*,\cdots,
\hat{u}^*)$ given in (\ref{eq27a}) is an open-loop $(\varepsilon_1,\varepsilon_2)$-Stackelberg equilibrium, 
where $\varepsilon_i=O(1/\sqrt{N})$, $i=1,2$.
\end{theorem}

\emph{Proof.} See Appendix A. \hfill{$\Box$}

\begin{theorem}\label{thm3.5}
For {Problem (PO)}, assume (A1)-(A4) hold, and $\xi_i,i=1,\cdots,N$ have the same variance. Under the control (\ref{eq27a}), the corresponding social cost 
is given by
\begin{equation}\label{eq426}
	J^{(N)}_{\rm soc}(u^*,u_0^*)=\mathbb{E}[|\xi_i|^2_{P(0)}+|\bar{\xi}|^2_{\bar{P}(0)}+|\bar{\xi}_0|^2_{K(0)}+2\bar{\xi}_0^TP_0(0)\bar{\xi}+2\varphi^T(0)\bar{\xi}+2\varphi^T_0(0)\bar{\xi}_0]+s_T,
\end{equation}
and the asymptotic cost of the leader is 
$
\lim_{N\to\infty}J_{0}(u_0^*,u^*)=\mathbb{E}\big[{\xi}^T_0y_0(0)+\bar{\xi}^T\bar{y}(0)\big] $, where
\begin{align}
	s_T
	=& \mathbb{E}\int_0^T\big[ (\bar{x}^T\bar{G}_0^TKD_0+\zeta_0^TD_0+\varphi_0^TB_0)u_0^*+u_0^*D_0KD_0u_0^*	\big]dt.
\end{align}

\end{theorem}
\emph{Proof.} See Appendix B. \hfill{$\Box$}

\section{Feedback Solutions to MF Leader-Follower Games}
\def\theequation{4.\arabic{equation}}
\setcounter{equation}{0}
In this section, we consider the feedback solution to the MF Stackelberg game (\ref{eq1})-(\ref{eq4}).
For simplicity, we consider the case that
$Q\geq 0,\ Q_0\geq 0,\ R>0$, $R_0>0$,  $H\geq 0$ and $H_0\geq 0$.

\subsection{The MF Social Control Problem for $N$ Followers}

Note that the leader plays against all  followers. Assume that the leader admits a feedback control of the following form
\begin{equation}
\label{eq35}
u_0=P_0x_0+\bar{P}x^{(N)},
\end{equation}
where $P_0$ and $\bar{P}$ are fixed. 
Thus, we have the following social control problem for $N$ followers.

\textbf{(P3)}: minimize $ {J}^{(N)}_{\rm soc}(u)$ over $u\in{\mathcal U}_{c}$, where $u_0=P_0x_0+\bar{P}x^{(N)}$ and
{\begin{equation}\label{eq7-a}
	{J}^{(N)}_{\rm soc}(u)=\frac{1}{N}\sum_{i=1}^N
	\mathbb{E}\int_0^{T}\!\!
	\Big\{\big|x_i
	-\Gamma x^{(N)}\!-\Gamma_1x_0\big|^2_{Q}+|u_i|^2_{R}\Big\}dt+\frac{1}{N}\sum_{i=1}^N\mathbb{E}\big[|x_i(T)-\hat{\Gamma}x^{(N)}(T)-\hat{\Gamma}_1x_0(T)|^2_{H}\big].
	\end{equation}}
	By examining the social cost variation, we obtain the optimal control laws for $N$ followers. The proof is similar to that of Theorem \ref{thm3.1}, and so omitted.
	\begin{theorem}\label{thm5.1}
Suppose $Q\geq0$, $R>0$ and $H\geq 0$. Assume the leader has the feedback control (\ref{eq35}). Then Problem (P3) has an optimal control in ${\mathcal U}_{c}$ if and only if  the following system of FBSDEs admits a set of adapted solutions $\{x_i,p_i,q_i^j,i,j=0,1,\cdots,N\}$:
\begin{equation}\label{eq59}
	\left\{
	\begin{aligned}
		dx_0=&\big[A_0x_0+B_0(P_0x_0+\bar{P}x^{(N)})+G_0x^{(N)}\big]dt+\big[C_0x_0+D_0(P_0x_0+\bar{P}x^{(N)})+\bar{G}_0x^{(N)}\big]dW_0,\cr
		dx_i=& (Ax_i+B\breve{u}_i+Gx^{(N)}+Fx_0)dt+(Cx_i+D\breve{u}_i+\bar{G}x^{(N)}+\bar{F}x_0)dW_i,\cr
		dp_0=&-\big[(A_0+B_0P_0)^Tp_0+F^Tp^{(N)}+(C_0+D_0P_0)^Tq_0^0+\bar{F}^Tq^{(N)}\cr
		&-Q_{\Gamma_1}^T x^{(N)}+\Gamma_1^TQ\Gamma_1 x_0)\big]
		+\sum_{j=0}^Nq_0^jdW_j,\\
		dp_i=&-\big[A^Tp_i+G^Tp^{(N)}+(G_0+ B_0\bar{P}) ^Tp_0+C^Tq_i^i+\bar{G}^Tq^{(N)}+(\bar{G}_0+D_0\bar{P})^Tq_0^0\cr
		&+Qx_i-Q_{\Gamma}x^{(N)}-Q_{\Gamma_1}\Gamma_1x_0\big]dt+\sum_{j=0}^Nq_i^jdW_j,\\
		x_0(0)&=\xi_{0},\ x_i(0)=\xi_{i},\ p_0(T)=-H_{\hat{\Gamma}_1}^T{x}^{(N)}(T)
		+\hat{\Gamma}_1^TH\hat{\Gamma}_1{x}_0(T),\cr
		p_i(T)&=H{x}_i(T)-H_{\hat{\Gamma}}{x}^{(N)}(T)-H_{\hat{\Gamma}_1} {x}_0(T),\ i=1,\cdots,N.
	\end{aligned}\right.
\end{equation}
Furthermore, the optimal controls of followers are given by
\begin{equation}\label{eq59a}
	\breve{u}_i=-R^{-1}(B^Tp_i+D^Tq_i^i),\ i=1,\cdots,N.
\end{equation}

\end{theorem}
\emph{Proof.}  See Appendix C. \hfill{$\Box$}

\begin{remark}
For the feedback solution case, the 
term $x^{(N)}$ appears in leader's dynamics. Distinct from the open-loop case, an additional costate $p_0$ is needed.
Indeed, as 
$u_i$ is perturbed with $\delta u_i$, the changing magnitude of $x^{(N)}$ is $O(\|\delta u_i\|/N)$, which causes the perturbation $O(\|\delta u_i\|)$ of $J_{\rm soc}(u)$.
This is evidently different from the game problem.
\end{remark}

Define
\begin{equation}\label{eq71}
\left\{\begin{aligned}
	&\dot{M}_N+A^TM_N+M_N^TA+C^TM_NC+Q-(B^TM_N+D^T\check{M}_NC)^T\Upsilon_N^{-1}\cr
	&\quad\times(B^TM_N+D^T\check{M}_NC)=0,\ M_N(T)=H,\cr
	&\dot{\bar{M}}_N+(A+G)^T\bar{M}_N+\bar{M}_N(A+G)+G^TM_N+M_NG+C^T\check{M}_N\bar{G}\cr
	&\quad+\bar{G}^T\check{M}_N(C+\bar{G})-Q_{\Gamma}+\bar{P}^TD_0^T\check{\Lambda}_N^0D_0\bar{P}+M_N^0(G_0+B_0\bar{P})+(G_0+B_0\bar{P})^T\bar{\Lambda}_N\cr
	&\quad-\big(B^TM_N+D^T\check{M}_NC\big)^T\Upsilon_N^{-1}(B^T\bar{M}_N+D^T\check{M}_N\bar{G})\cr
	&\quad-\big(B\bar{M}_N+D^T\check{M}_N\bar{G}\big)^T\Upsilon_N^{-1}\big(B^TM_N+D^T\check{M}_NC\big)\cr
	&\quad-\big(B\bar{M}_N+D^T\check{M}_N\bar{G}\big)^T\Upsilon_N^{-1}\big(B^T\bar{M}_N+D^T\check{M}_N\bar{G}\big)=0,\bar{M}_N(T)=-H_{\hat{\Gamma}},\cr
	&\dot{M}^0_N+(A+G)^TM_N^0+M^0_N(A_0+B_0P_0)+(M_N+\bar{M}_N)F+(G_0+B_0\bar{P})^T\Lambda_N^0\cr
	&\quad-[B^T(M_N+\bar{M}_N)+D^T\check{M}_N(C+\bar{G})]^T\Upsilon^{-1}_N(B^TM^0_N+D^T\check{M}\bar{F})\cr
	&\quad+(C+\bar{G})^T\check{M}_N\bar{F}+(\bar{G}_0+D_0\bar{P})^T\check{\Lambda}_N^0(C_0+D_0P_0)+(\Gamma-I)^TQ\Gamma_1=0,\cr
	&\quad M^0_N(T)=(\hat{\Gamma}-I)^TH\hat{\Gamma}_1,
\end{aligned}\right.
\end{equation}
\begin{equation}\label{eq72}
\left\{\begin{aligned}
	&\dot{\Lambda}^0_N+\Lambda^0_N(A_0+B_0P_0)+(A_0+B_0P_0)^T\Lambda^0_N+(C_0+D_0P_0)^T\check{\Lambda}_N^0(C_0+D_0P_0)\cr
	&\quad-(B^T\bar{\Lambda}^T_N+D^T\check{M}_N\bar{F})^T\Upsilon_N^{-1}(B^TM^0_N+D^T\check{M}_N\bar{F})\cr
	&\quad+\bar{\Lambda}_NF+F^TM_N^0+\bar{F}^T\check{M}\bar{F}+\Gamma_1^TQ\Gamma_1=0, \ \Lambda^0_N(T)=\hat{\Gamma}_1^TH\hat{\Gamma}_1,\cr
	&\dot{\bar{\Lambda}}_N+\bar{\Lambda}_N(A+G)+(A_0+B_0P_0)^T\bar{\Lambda}_N+F^T(M_N+\bar{M}_N)+\Lambda^0_N(G_0+B_0\bar{P})\cr
	&\quad-(B^T\bar{\Lambda}_N^T+D^T\check{M}_N\bar{F})^T\Upsilon^{-1}_N[B^T(M_N+\bar{M}_N)+D^T\check{M}_N(C+\bar{G})]\cr
	&\quad+\bar{F}^T\check{M}_N(C+\bar{G})+\Gamma_1^T Q(\Gamma-I)=0, \ \bar{\Lambda}_N(T)= \hat{\Gamma}_1^TH(\hat{\Gamma}-I).
\end{aligned}\right.
\end{equation}

\begin{proposition}\label{prop5.1}
Assume (A1) holds, and (\ref{eq71})-(\ref{eq72}) admit solutions, respectively. Then, Problem (P3) admits a feedback solution (\ref{eq67}).
\end{proposition}
{\emph{Proof.} Let $p_0=\Lambda^0_Nx_0+\bar{\Lambda}_Nx^{(N)}$, and  $p_i=M_N x_i+\bar{M}_Nx^{(N)}+M^0_Nx_0,\ i=1,\cdots,N.$
Denote  $\breve{u}^{(N)}=\frac{1}{N}\sum_{i=1}^N\breve{u}_i$. By applying It\^{o}'s formula to $p_i$, we have
\begin{align}\label{eq66}
dp_i=&\dot{M}_Nx_idt+M_N\big[(Ax_i+B\check{u}_i+Gx^{(N)}+Fx_0)dt+(Cx_i+D\check{u}_i+\bar{G}x^{(N)}+\bar{F}x_0)dW_i\big]\\
&+\dot{\bar{M}}_Nx^{(N)}+\bar{M}_N\Big[(A+ G)x^{(N)}+B\check{u}^{(N)}+Fx_0)dt+\frac{1}{N}\sum_{j=1}^N(Cx_j+D\check{u}_j+\bar{G}x^{(N)}+\bar{F}x_0)dW_j\Big]\cr
&+\dot{M}_N^0x_0dt+M_N^0\big[\big((A_0+B_0P_0)x_0+(G_0+B_0\bar{P})x^{(N)}\big)dt+\big((C_0+D_0P_0)x_0+(\bar{G}_0+D_0\bar{P})x^{(N)}\big)dW_0\big]\cr
=&-\big[A^T(M_N x_i+\bar{M}_Nx^{(N)}+M_N^0x_0)+G^T\big((M_N+\bar{M}_N)x^{(N)}+M_N^0x_0\big)+(G_0+B_0\bar{P})^Tp_0\cr
&+C^Tq_i^i+\bar{G}^Tq^{(N)}+(\bar{G}_0+D_0\bar{P})^Tq_0^0+Qx_i-Q_{\Gamma}x^{(N)}+(\Gamma-I)^TQ\Gamma_1x_0\big]dt+\sum_{j=0}^Nq_i^jdW_j,\nonumber
\end{align}
which together with (\ref{eq59}) implies
\begin{equation}
\label{eq67a}
\begin{aligned}
	q_i^i=& \big(M_N+\frac{1}{N}\bar{M}_N\big)(Cx_i+D\check{u}_i+\bar{G}x^{(N)}+\bar{F}x_0),\cr
	q_i^j=&\frac{1}{N}\bar{M}_N(Cx_j+D\check{u}_j+\bar{G}x^{(N)}+\bar{F}x_0),\ j\not=i.
\end{aligned}
\end{equation}
By (\ref{eq59a}), we have for any $i=1,\cdots,N$,
\begin{align*}
R\breve{u}_i+B^T(M_N x_i+\bar{M}_Nx^{(N)}+M_N^0x_0)+D^T\big(M_N+\frac{1}{N}\bar{M}_N\big)(Cx_i+D\check{u}_i+\bar{G}x^{(N)}+\bar{F}x_0)=0.	
\end{align*}
This leads to
\begin{equation}\label{eq67}
\breve{u}_i=-\Upsilon_N^{-1}\big[(B^TM_N+D^T\check{M}_NC)x_i+(B^T\bar{M}_N+D^T\check{M}_N\bar{G})x^{(N)}+(B^TM_N^0+D^T\check{M}_N\bar{F})x_0\big],
\end{equation}
where $\check{M}_N\stackrel{\Delta}{=}M+\frac{1}{N}\bar{M}_N$ and $\Upsilon_N\stackrel{\Delta}{=}R+D^T\check{M}_ND$.
Denote $\check{\Lambda}_N^0\stackrel{\Delta}{=}{\Lambda}_N^0+\frac{1}{N}\bar{\Lambda}_N$.
Applying It\^{o}'s formula to $p_0$, we obtain
\begin{align}\label{eq68a}
dp_0
\!=&\dot{\Lambda}_N^0x_0dt\!+\!\Lambda^0_N\big[\big(A_0x_0\!+\!B_0(P_0x_0\!+\!\bar{P}x^{(\!N\!)})\!+\!{G}_0x^{(\!N\!)}\big)dt
\!+\!\big(C_0x_0\!+\!D_0(P_0x_0\!+\!\bar{P}x^{(\!N\!)})\!+\!\bar{G}_0x^{(\!N\!)}\big)dW_0\big]\cr
&+\!\dot{\bar{\Lambda}}_Nx^{(\!N\!)}\!+\!\bar{\Lambda}_N\Big[\big((A\!+\!G)x^{(\!N\!)}\!+\!B\check{u}^{(\!N\!)}\!+\!Fx_0\big)dt
\!+\!\frac{1}{N}\sum_{j=1}^N(Cx_j\!+\!D\check{u}_j\!+\!\bar{G}x^{(\!N\!)}\!+\!\bar{F}x_0)dW_j\Big],
\end{align}
which together with (\ref{eq59}) implies
\begin{equation}
\label{eq68}
\begin{aligned}
	&q_0^0=\check{\Lambda}_N^0\big(C_0x_0+D_0(P_0x_0+\bar{P}x^{(N)})+\bar{G}_0x^{(N)}\big),\cr &q_0^j=\frac{1}{N}\bar{\Lambda}\big(C_0x_0+D_0(P_0x_0+\bar{P}x^{(N)})+\bar{G}_0x^{(N)}\big),\ j>0.
\end{aligned}
\end{equation}
Applying
(\ref{eq67a}), (\ref{eq67}) and (\ref{eq68}) into (\ref{eq66}), we obtain (\ref{eq71}). Applying
(\ref{eq67a}), (\ref{eq67}) and (\ref{eq68}) into (\ref{eq68a}), we have (\ref{eq72}). Based on Theorem \ref{thm5.1} and the above discussion,  the proposition follows. \hfill{$\Box$}
}

\begin{remark}
Note that the social problem (P3) is essentially an optimal control problem. The  feedback solution to
Problem (P3) is equivalent to the feedback representation of its open-loop solution.
\end{remark}

We now introduce the following set of equations:
\begin{equation}\label{eq73}
\left\{\begin{aligned}
	&\dot{M}+A^TM+M^TA+C^TMC+Q-(B^TM+D^T{M}C)^T\Upsilon^{-1}\cr
	&\quad \times (B^TM+D^T{M}C)=0, \ M(T)=H,\cr
	&\dot{\bar{M}}+(A+G)^T\bar{M}+\bar{M}(A+G)+G^TM+MG+C^T{M}\bar{G}+\bar{G}^T{M}(C+\bar{G})\cr
	&\quad-\big(B^TM+D^T{M}C\big)^T\Upsilon^{-1}(B^T\bar{M}+D^T{M}\bar{G})+\bar{P}^TD_0^T{\Lambda}^0D_0\bar{P}+(G_0+B_0\bar{P})^T\bar{\Lambda}\cr
	&\quad-\big(B\bar{M}+D^T{M}\bar{G}\big)^T\Upsilon^{-1}\big(B^TM+D^T{M}C\big)-Q_{\Gamma}+M^0(G_0+B_0\bar{P})\cr
	&\quad-\big(B\bar{M}+D^T{M}\bar{G}\big)^T\Upsilon^{-1}\big(B^T\bar{M}+D^T{M}\bar{G}\big)=0,\ \bar{M}(T)=-H_{\hat{\Gamma}},\cr
	&\dot{M}^0+(A+G)^TM^0+M^0(A_0+B_0P_0)+(M+\bar{M})F+(G_0+B_0\bar{P})^T\Lambda^0\cr
	&\quad-[B^T(M+\bar{M})+D^T{M}(C+\bar{G})]^T\Upsilon^{-1}(B^TM^0+D^T{M}\bar{F})+(C+\bar{G})^T{M}\bar{F}\cr
	&\quad+(\bar{G}_0+D_0\bar{P})^T{\Lambda}^0(C_0+D_0P_0))+(\Gamma-I)^TQ\Gamma_1=0,\ M^0(T)=(\hat{\Gamma}-I)^TH\hat{\Gamma}_1,\cr
	&\dot{\Lambda}^0+\Lambda^0(A_0+B_0P_0)+(A_0+B_0P_0)^T\Lambda^0+(C_0+D_0P_0)^T{\Lambda}^0(C_0+D_0P_0)\cr
	&\quad-(B^T\bar{\Lambda}^T+D^T{M}\bar{F})^T\Upsilon^{-1}(B^TM^0+D^T{M}\bar{F})+\bar{\Lambda}F+F^TM^0+\bar{F}^T{M}\bar{F}\cr
	&\quad+\Gamma_1^TQ\Gamma_1=0, \ \Lambda^0(T)=\hat{\Gamma}_1^TH\hat{\Gamma}_1,\cr
	&\dot{\bar{\Lambda}}+\bar{\Lambda}(A+G)+(A_0+B_0P_0)^T\bar{\Lambda}+F^T(M+\bar{M})+\Lambda^0(G_0+B_0\bar{P})\cr
	&\quad-(B^T\bar{\Lambda}^T+D^T{M}\bar{F})^T\Upsilon^{-1}[B^T(M+\bar{M})+D^T{M}(C+\bar{G})]+\bar{F}^T{M}(C+\bar{G})\cr
	&\quad+\Gamma_1^T Q(\Gamma-I)=0, \ \bar{\Lambda}(T)= \hat{\Gamma}_1^TH(\hat{\Gamma}-I),
\end{aligned}\right.
\end{equation}
where $\Upsilon
\stackrel{\Delta}{=}R+D^TMD$. From observation, we find that $M,\bar{M},\Lambda^0$ are symmetric and $M^0=\bar{\Lambda}^T$.  For further analysis, we assume

\textbf{(A5)}  (\ref{eq73}) admits a  solution $(M,\bar{M},M^0,\Lambda^0,\bar{\Lambda})$.

\begin{remark}\label{rem5.1}
If (A5) holds, then by the continuous dependence of solutions on the parameter (see e.g. \cite[Theorem 3.5]{K02} or
\cite[Theorem 4]{HZ20}),
we obtain that for sufficiently large $N$, (\ref{eq71}) and (\ref{eq72}) admit solutions, respectively.
\end{remark}

After applying the strategies of followers (\ref{eq67}), we have
\begin{align}
d{x}_i
=&\big[(A-B\Upsilon_N^{-1}\Psi_N){x}_i+(G-B\Upsilon_N^{-1}\bar{\Psi}_N){x}^{(N)}+(F-B\Upsilon_N^{-1}\Psi_N^0)x_0\big]dt\\
&+\big[(C-D\Upsilon_N^{-1}\Psi_N){x}_i+(\bar{G}-D\Upsilon_N^{-1}\bar{\Psi}_N){x}^{(N)}+
({\bar{F}}-D\Upsilon_N^{-1}\Psi_N^0)x_0\big]dW_i,\nonumber
\end{align}
where $\Psi_N\stackrel{\Delta}{=}B^{T}M_N+D^T\check{M}_NC$, $\bar{\Psi}_N=B^{T}\bar{M}_N+D^T\check{M}_N\bar{G}$, and $\Psi_N^0=B^T{M}_N^0+D^T\check{M}_N\bar{F}$.
This leads to
\begin{align*}
d{x}^{(N)}=&\big[\big(A+G-B\Upsilon_N^{-1}(\Psi_N+\bar{\Psi}_N)\big){x}^{(N)}+(F-B\Upsilon_N^{-1}\Psi_N^0)x_0\big]dt\cr
&+\frac{1}{N}\sum_{i=1}^N\big[(C-D\Upsilon_N^{-1}\Psi_N){x}_i+(\bar{G}-D\Upsilon_N^{-1}\bar{\Psi}_N){x}^{(N)}+
({\bar{F}}-D\Upsilon_N^{-1}\Psi_N^0)x_0\big]dW_i.
\end{align*}
For a sufficiently large $N$, by Remark \ref{rem5.1} and the law of large numbers,  ${x}^{(N)}$ can be  approximated by the MF function $\bar{x}$, which
satisfies
\begin{align}\label{eq75a}
d\bar{x}=&\big[\big(A+G-B\Upsilon^{-1}(\Psi+\bar{\Psi})\big)\bar{x}+
(F-B\Upsilon^{-1}\Psi^0)x_0\big]dt,
\end{align}
with
\begin{equation}\label{eq518}
\begin{aligned}
	&\Psi\stackrel{\Delta}{=}B^{T}M+D^T{M}C,
	\ \bar{\Psi}\stackrel{\Delta}{=}B^{T}\bar{M}+D^T{M}\bar{G},\cr
	&\Psi^0\stackrel{\Delta}{=}B^T{M}^0+D^T{M}\bar{F}.
\end{aligned}
\end{equation}
Based on Proposition \ref{prop5.1}, one can construct the decentralized feedback strategies for followers:
\begin{equation}\label{eq75}
\hat{u}_i=-\Upsilon^{-1}(\Psi x_i+\bar{\Psi}\bar{x}+\Psi^0x_0).
\end{equation}
\subsection{Optimization for the Leader}

After applying the strategies (\ref{eq75}) of followers,
we have the optimal control problem for the leader.

\textbf{(P4)}: minimize $ {J}_{0}(u_0,\hat{u}(u_0))$ over $u_0\in \mathcal{U}_0$,
where
\begin{align*}
&J_0(u_0,\hat{u}(u_0))=\mathbb{E}\int_{0}^{T}\big[|x_0- \Gamma_0\hat{x}^{(N)}|_{Q_0}^{2}+|u_0|_{R_0}^{2}\big]dt+\mathbb{E}\big[|x_0(T)-\hat{\Gamma}_0x^{(N)}(T)|^2_{H_0}\big],\cr
&dx_0=\big(A_0x_0+B_0u_0+G_0\hat{x}^{(N)}\big)dt
+\big(C_0 x_0+D_0u_0+\bar{G}_0\hat{x}^{(N)}\big)dW_0,\ x_0(0)=\xi_0, \label{eq50a}
\\ 
&d\hat{x}_i=\big[(A-B\Upsilon^{-1}\Psi)\hat{x}_i+G\hat{x}^{(N)} -B\Upsilon^{-1}\bar{\Psi}\bar{x}+(F-B\Upsilon^{-1}\Psi^0)x_0\big]dt\cr
&\hspace{2em}+\big[(C-D\Upsilon^{-1}\Psi)\hat{x}_i+\bar{G}\hat{x}^{(N)}-D\Upsilon^{-1}\bar{\Psi}\bar{x}+(\bar{F}-D\Upsilon^{-1}\Psi^0)x_0\big]dW_i,\ \hat{x}_i(0)=\xi_i.
\end{align*}
Since $\{W_i(t)\}$ and $\{x_i(0)\}$ are independent sequences, for a sufficiently large $N$, it is plausible to replace $\hat{x}^{(N)}$ by $\bar{x}$, which evolves from (\ref{eq75a}). 
In view of (\ref{eq35}), suppose that the decentralized feedback solution for the leader has the following form
$u_0(t)=P_0(t)x_0+\bar{P}(t)\bar{x}, \ 0\leq t\leq T.$
Then, we have the following optimal control problem for the leader.

\textbf{(P4$^{\prime}$)}: minimize $ \bar{J}_{0}(P_0,\bar{P})$ over $P_0,\bar{P}\in C(0,T;\mathbb{R}^{m\times n})$, where
\begin{equation*}
\left\{\begin{aligned}
	&\bar{J}_0(P_0,\bar{P})=\mathbb{E}\int_{0}^{T}\big[|x_0- \Gamma_0\bar{x}|_{Q_0}^{2}+|P_0x_0+\bar{P}\bar{x}|_{R_0}^{2}\big]dt+\mathbb{E}\big[|x_0(T)-\hat{\Gamma}_0\bar{x}(T)|^2_{H_0}\big],
	\cr
	&d\bar{x}_0=\big[(A_0+B_0P_0)\bar{x}_0+(G_0+B_0\bar{P})\bar{x}\big]dt+\big[\big(C_0+D_0P_0\big)\bar{x}_0+(\bar{G}_0+D_0\bar{P})\bar{x}\big]dW_0,\ \bar{x}_0(0)=\xi_0,  
	\\
	&d\bar{x}=\big [\big(A+G-B\Upsilon^{-1}(\Psi+\bar{\Psi})\big)\bar{x}+(F-B\Upsilon^{-1}\Psi^0)\bar{x}_0\big]dt, \ \bar{x}(0)=\bar{\xi}.
\end{aligned}\right.
\end{equation*}
Let $X_0=\mathbb{E}[\bar{x}_0\bar{x}_0^T]$, $\bar{X}=\mathbb{E}[\bar{x}\bar{x}^T]$ and $Y=\mathbb{E}[\bar{x}_0\bar{x}^T]$.
Then, by It\^{o}'s formula \cite{YZ99}, we obtain
\begin{align}
\frac{d\bar{X}_0}{dt}=&(A_0+B_0P_0)X_0+X_0(A_0+B_0P_0)^T+(G_0+B_0\bar{P})Y^T+Y(G_0+B_0\bar{P})^T\\
&+(C_0+D_0P_0)X_0(C_0+D_0P_0)^T+(C_0+D_0P_0)Y(\bar{G}_0+D_0\bar{P})^T\cr
&+(\bar{G}_0+D_0\bar{P})Y^T(C_0+D_0P_0)^T+(\bar{G}_0+D_0\bar{P})\bar{X}(\bar{G}_0+D_0\bar{P})^T,\cr
\frac{d\bar{X}}{dt}=&(A+G-B\Upsilon^{-1}(\Psi+\bar{\Psi}))\bar{X}+\bar{X}(A+G-B\Upsilon^{-1}(\Psi+\bar{\Psi}))^T\\
&+(F-B\Upsilon^{-1}\Psi^0)Y+Y^T(F-B\Upsilon^{-1}\Psi^0)^T,\cr
\frac{dY}{dt} =&Y(A+G-B\Upsilon^{-1}(\Psi+\bar{\Psi}))^T+\bar{X}_0(F-B\Upsilon^{-1}\Psi_0)^T\cr
&+(A_0+B_0P_0)Y+(G_0+B_0\bar{P})\bar{X}.
\end{align}
Meanwhile, the cost function of the leader can be rewritten as
\begin{align*}
\bar{J}_0(P_0,\bar{P})=&\int_0^Ttr\big(Q_0X_0-Q_0\Gamma_0Y^T-\Gamma_0^TQ_0Y+\Gamma_0^TQ_0\Gamma_0\bar{X}\cr
&+P_0^TR_0P_0X_0+\bar{P}^TR_0P_0Y+P_0^TR_0\bar{P}Y^T+\bar{P}^TR_0\bar{P}\bar{X}\big)dt\cr
&+tr \big[H_0{X}_0(T)-H_0\hat{\Gamma}_0Y^T(T)-\hat{\Gamma}_0^TH_0Y(T)+\hat{\Gamma}_0^TH_0\hat{\Gamma}_0\bar{X}(T)\big].
\end{align*}
Denote
$\hat{A}_0\stackrel{\Delta}{=}A_0+B_0P_0,\ \hat{C}_0\stackrel{\Delta}{=}C_0
+D_0P_0,\ \hat{F}\stackrel{\Delta}{=}F-B\Upsilon^{-1}\Psi^0,\  \hat{A}\stackrel{\Delta}{=}A+G-B\Upsilon^{-1}(\Psi+\bar{\Psi}).$
Define the Hamiltonian function of the leader as follow:
\begin{align*}
&H(P_0,\bar{P},\Theta_1,\Theta_2,\Theta_3)\cr
=&tr\Big(Q_0X_0-Q_0\Gamma_0Y^T-\Gamma_0^TQ_0Y+\Gamma_0^TQ_0\Gamma_0\bar{X}+P_0^TR_0P_0X_0+\bar{P}^TR_0P_0Y+P_0^TR_0\bar{P}Y^T\cr
&+\bar{P}^TR_0\bar{P}\bar{X}+[\hat{A}_0X_0+X_0\hat{A}_0^T+(G_0+B_0\bar{P})Y^T+Y(G_0+B_0\bar{P})^T+\hat{C}_0Y(\bar{G}_0+D_0\bar{P})^T\cr
&+\hat{C}_0X_0\hat{C}_0^T+(\bar{G}_0+D_0\bar{P})Y^T\hat{C}_0^T+(\bar{G}_0+D_0\bar{P})\bar{X}(\bar{G}_0+D_0\bar{P})^T]\Theta_1^T+[\hat{A}\bar{X}+\bar{X}\hat{A}^T+\hat{F}Y+Y^T\hat{F}^T]\Theta_2^T\cr
&+\big[Y\hat{A}^T+X_0\hat{F}^T+\hat{A}_0Y+(G_0+B_0\bar{P})\bar{X}\big]\Theta_3^T+\big[Y\hat{A}^T+X_0\hat{F}^T+\hat{A}_0Y+(G_0+B_0\bar{P})\bar{X}\big]^T\Theta_3\Big).
\end{align*}
By the matrix maximum principle \cite{A68}, we obtain the following adjoint equations:
\begin{equation}\label{eq80}
\left\{ \begin{aligned}
	-\dot{\Theta}_1=\frac{\partial H}{\partial X_0}=&Q_0+P_0^TR_0P_0+\hat{A}_0^T\Theta_1+\Theta_1\hat{A}_0+\hat{C}_0^T\Theta_1\hat{C}_0+\hat{F}^T\Theta_3^T+\Theta_3\hat{F},\\
	- \dot{\Theta}_2=\frac{\partial H}{\partial \bar{X}}=&\Gamma_0^TQ\Gamma_0+\bar{P}^TR_0\bar{P}+\hat{A}^T\Theta_2+\Theta_2\hat{A}+\Theta_3^T(G_0+B_0\bar{P})+(G_0+B_0\bar{P})^T\Theta_3\cr
&+(\bar{G}_0+D_0\bar{P})^T\Theta_1(\bar{G}_0+D_0\bar{P}),\\
	-  \dot{\Theta}_3=\frac{\partial H}{\partial Y}=&P_0^TR_0\bar{P}-Q_0\Gamma_0+\Theta_1(G_0+B_0\bar{P})+\hat{F}^T\Theta_2+\hat{C}_0^T\Theta_1(\bar{G}_0+D_0\bar{P})+\Theta_3\hat{A}+\hat{A}_0^T\Theta_3,&
\end{aligned}\right.
\end{equation}
with the stationarity conditions
\begin{align}\label{eq831}
0=&\frac{\partial H}{\partial P_0}=2\big[R_0P_0X_0+R_0\bar{P}Y^T+B_0^T\Theta_1X_0+D_0^T\Theta_1\hat{C}_0X_0+D_0^T\Theta_1(\bar{G}_0+D_0\bar{P})Y^T +B_0^T\Theta_3Y^T\big],\\
\label{eq84}
0=&\frac{\partial H}{\partial \bar{P}}= 2\big[R_0P_0Y+R_0\bar{P}\bar{X}+B_0^T\Theta_1Y^T+D_0^T\Theta_1\hat{C}_0Y^T+D_0^T\Theta_1(\bar{G}_0+D_0\bar{P})\bar{X}+B_0^T\Theta_3\bar{X}\big].
\end{align}
Note that $\Theta_1$ and $\Theta_2$ are symmetric matrices. From  (\ref{eq831}) and (\ref{eq84}), we obtain
$\left\{ \begin{aligned}
	P_0&=-\Upsilon_0^{-1}\Psi_4\\
	\bar{P}&=-\Upsilon_0^{-1}\Psi_5,
\end{aligned}\right.$
where
\begin{equation}\label{eq433}
  \Upsilon_0{=}R_0+D_0^T\Theta_1D_0, \ \Psi_4=B_0^T\Theta_1+D_0^T\Theta_1C_0,\ \Psi_5=B_0^T\Theta_3+D_0^T\Theta_1\bar{G}_0.
  \end{equation}
Applying this into (\ref{eq80}), we have
\begin{equation}\label{eq832}
\left\{ \begin{aligned}
	&\dot{\Theta}_1+A_0^T\Theta_1+\Theta_1A_0+{C}_0^T\Theta_1{C}_0-\Psi_4^T\Upsilon_0^{-1}\Psi_4+\hat{F}^T\Theta_3^T+\Theta_3\hat{F}+Q_0=0,\ \Theta_1(T)=H_0,\cr
	&\dot{\Theta}_2+\hat{A}^T\Theta_2+\Theta_2\hat{A}-\Psi_5^T\Upsilon_0^{-1}\Psi_5 +\Gamma_0^TQ\Gamma_0+\bar{G}_0^T\Theta_1\bar{G}_0+\Theta_3^TG_0+G_0^T\Theta_3=0, \ \Theta_2(T)=\hat{\Gamma}_0^TH_0\hat{\Gamma}_0,\\
	&\dot{\Theta}_3+\Theta_3\hat{A}+A_0^T\Theta_3-\Psi_4^T\Upsilon_0^{-1}\Psi_5 +\hat{F}^T\Theta_2-Q_0\Gamma_0+\Theta_1G_0+C_0^T\Theta_1\bar{G}_0=0,\ \Theta_3(T)=-H_0\hat{\Gamma}_0.
\end{aligned}\right.
\end{equation}

Based on the above discussions, we may construct the following feedback strategies:
\begin{equation}\label{eq33a}
\left\{\begin{aligned}
	\hat{u}_0=&-\Upsilon_0^{-1}\big[\Psi_4x_0+\Psi_5\bar{x}\big],\cr
	\hat{u}_i=&-\Upsilon^{-1}(\Psi x_i+\bar{\Psi}\bar{x}+\Psi^0x_0), \ i=1,\cdots,N,
\end{aligned}\right.
\end{equation}
where $\Upsilon_0, \Psi_4,\Psi_5$ are given by
(\ref{eq433}),
and $\Psi$, $\bar{\Psi}$, $\Psi^0$ are given by (\ref{eq518}).

%
%
%

\begin{theorem}\label{thm4.4}
For {Problem (PF)}, assume (A1) holds; (\ref{eq73}) and (\ref{eq832}) admit a set of solutions. Then, the strategy (\ref{eq33a}) is a feedback ($\epsilon_1,\epsilon_2$)-Stackelberg equilibrium, where $\epsilon_1=\epsilon_2=O(\frac{1}{\sqrt{N}})$. Furthermore, assume that $\xi_i,i=1,\cdots,N$ have the same variance. Then, the asymptotic average social cost of followers is given by
\begin{align*}
	\lim_{N\to\infty}\frac{1}{N}J_{\rm soc}(\hat{u},\hat{u}_0)=\mathbb{E}[|\xi_i|^2_{M(0)}+|\bar{\xi}|^2_{\bar{M}(0)}+2\xi_0^T\bar{\Lambda}(0)\xi_i+|\xi_0|^2_{\Lambda_0(0)}],
\end{align*}
and
\begin{align*}
	\lim_{N\to\infty}J_{0}(\hat{u},\hat{u}_0)=\mathbb{E}[\xi_0^T\Theta_1(0)\xi_0+\bar{\xi}^T\Theta_2(0)\bar{\xi}+\bar{\xi}^T\Theta_3(0)\xi_0].
\end{align*}
\end{theorem}

\emph{Proof. }
See Appendix C. \hfill{$\Box$}


\section{Simulation}

In this section, we give a numerical example to compare the performances of the open-loop and feedback solutions.
The simulation parameters are listed in Table \ref{tab21}.  

\begin{figure}[H]
\includegraphics[scale=0.7]{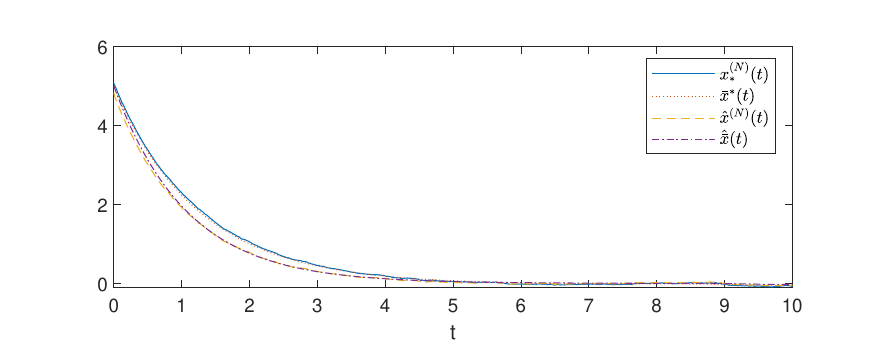}
\centering
\vspace{-1em}
\caption{  The solution to the Riccati equation \eqref{eq18}, and $P_{i,j}$ is the entry in $i$th row $j$th column of $\mathcal{P}$.}
\label{f1}
\end{figure}

\begin{figure}[H]
\includegraphics[scale=0.7]{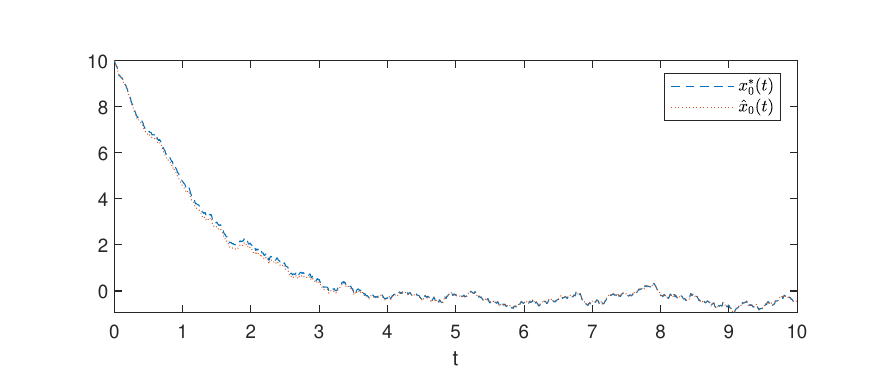}
\centering
\vspace{-1em}
\caption{ The solutions to \eqref{eq73}  and \eqref{eq832}.}
\label{f2}
\end{figure}

\begin{figure}[H]
\includegraphics[scale=0.7]{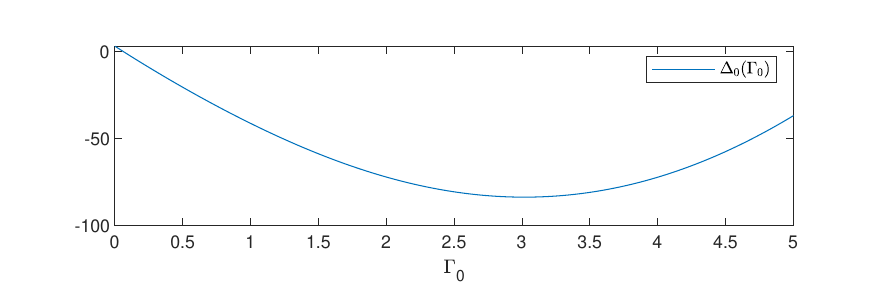}
\centering
\vspace{-1em}
\caption{  Followers' state averages and MF effects under  open-loop and feedback controls.}
\label{f3}
\end{figure}

\begin{figure}[H]
\includegraphics[scale=0.7]{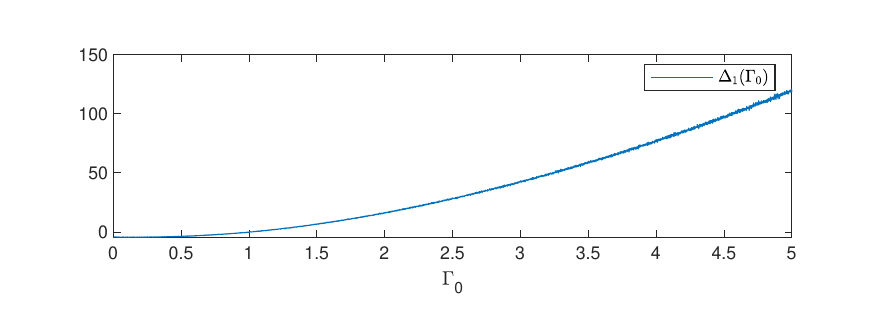}
\centering
\vspace{-1em}
\caption{ States of the leader under  open-loop and feedback controls.}
\label{f4}
\end{figure}

\vspace{-1em}
\begin{table}[H]
\centering
\setlength{\tabcolsep}{1.2mm}
\caption{Simulation parameters}
\label{tab21}       
\begin{tabular}{cccccccccccccccc}
	\hline\noalign{\smallskip}
	& &	$A_0$&$B_0$&$C_0$&$D_0$&$\Gamma_0$&$Q_0$&$R_0$&$\hat{\Gamma}_0$&$H_0$\\
	\noalign{\smallskip}\hline\noalign{\smallskip}
	& &	$-10$&$1$&$-0.5$&$0.5$&$1$&$1$&$1$&$1$&$2$ \\
	\noalign{\smallskip}\hline\noalign{\smallskip}
	$A$&$B$&$G$&$F$&$C$&$D$&$\bar{G}$&$\bar{F}$&$\Gamma$&$\Gamma_1$&$Q$&$R$&$\hat{\Gamma}$&$\hat{\Gamma}_1$&$H$\\
	\noalign{\smallskip}\hline\noalign{\smallskip}
	$-2$&$1$&$1$&$1$&$-0.2$&$0.2$&$0.2$&$0.2$&$1$&$1$&$1$&$1$&$1$&$1$&$2$\\
	\noalign{\smallskip}\hline\noalign{\smallskip}
\end{tabular}
\end{table}
\vspace{-1em}
Consider a multi-agent system with $1$ leader and $100$ followers. The initial distributions of states for the leader and followers satisfy  normal distributions $N(10,2)$ and $N(5,1)$, respectively. The decentralized open-loop control \eqref{eq27a} is given  by solving    
\eqref{eq9b}, \eqref{eq10b}, \eqref{eq7} and (\ref{eq18}). The solution to the Riccati equation \eqref{eq18} is shown in Fig. \ref{f1}.  The decentralized feedback strategy \eqref{eq33a} is obtained by solving \eqref{eq73} and \eqref{eq832}. The solutions to \eqref{eq73}  and \eqref{eq832} are shown in Fig. \ref{f2}.  Fig. \ref{f3} gives the curves of followers' state averages and MF effects under  open-loop and feedback solutions. Fig. \ref{f4} shows the state  trajectories of the leader under the two solutions. It can be seen  that  state averages approximate MF effects well under both solutions, and the state average under open-loop control is larger than the one under feedback control.
\vspace{-1em}
\vspace{-2em}
\vspace{-2em}
\vspace{-2em}

\section{Concluding Remarks}
This paper studies open-loop and feedback solutions of MF-LQG Stackelberg games with multiplicative noise. By
decoupling MF FBSDEs and applying MF approximations, we obtain a set of  open-loop
controls of players and a
set of decentralized feedback strategies, respectively. Furthermore, the corresponding optimal costs of all players are explicitly given in terms of the solutions to two Riccati equations, respectively.
{A challenge is computing the system of Riccati equations for feedback strategies. A possible approach is resorting to reinforcement learning 
even if dynamics are partially unknown.}

\appendix

\section{Proof of Theorems \ref{thm3.1} and \ref{thm4.2}}\label{app a}
\def\theequation{A.\arabic{equation}}
\setcounter{equation}{0}

\emph{Proof of Theorem \ref{thm3.1}.}
Suppose that $\{\check{u}_i,i=1,\cdots,N\}$ is a candidate of optimal control to Problem (P1).
		 Denote by $\check{x}_i$ the state of agent $i$ under the optimal control $\check{u}_i$. For any $u_i\in L^2_{{\mathcal F}}(0, T; \mathbb{R}^r) $ and $\theta\in \mathbb{R}\ (\theta \not= 0)$, let $u_i^{\theta}=\check{u}_i+\theta u_i$. Denote by $x_i^{\theta},  i=0,1,\cdots,N$ the solution of the following perturbed state equation:
	$$ \begin{aligned}
dx_0^{\theta}\!=&\big[A_0x_0^{\theta} +B_0\check{u}_0
+\frac{G_0}{N}\sum_{i=1}^Nx^{\theta}_i\big]dt
+\big[C_0x_0^{\theta} +D_0\check{u}_0
+\frac{\bar{G}_0}{N}\sum_{i=1}^Nx^{\theta}_i\big]dW_0, x_0^{\theta}(0)\!=\!\xi_{0},\cr
		  dx_i^{\theta}=&\big(Ax_i^{\theta}+B(\check{u}_i+\theta u_i)+\frac{1}{N}\sum_{i=1}^NGx^{\theta}_i+Fx_0\big)dt+\big(Cx_i^{\theta}+Du_i^{\theta}+\frac{1}{N}\sum_{i=1}^N\bar{G}x^{\theta}_i+\bar{F}x_0\big) dW_i,\ x_i^{\theta}(0)=\xi_{i}.
		 \end{aligned}$$ Let $y_i=(x_i^{\theta}-\check{x}_i)/\theta$. 
	It can be verified that
	 $z_i$ satisfies
\begin{equation*}
\left\{ \begin{aligned}
	dz_0=&\big[A_0z_0+G_0z^{(N)}\big]dt+\big[C_0z_0+\bar{G}_0z^{(N)}\big]dW_0,\  z_0(0)=0,\cr
	dz_i=&[A y_i+B u_i+Gy^{(N)}+Fy_0]dt+[C y_i+D u_i+\bar{G}y^{(N)}+\bar{F}y_0]dW_i,\  y_i(0)=0, \ i=1,2,\cdots,N.
\end{aligned}\right.
\end{equation*}
 We have
 \begin{equation}
   \label{eqA1}
J_{\rm  soc}(\check{u}+\theta u)-J_{\rm soc}(\check{u})=2\theta I_1+\theta^2 I_2,
 \end{equation} where
    \begin{align}\label{eq8a}
I_1=&\sum_{i=1}^N \mathbb{E}\int_0^T\big[\check{x}_i^TQy_i-(\check{x}^{(N)})^T{Q_{\Gamma}}y^{(N)}- x_0^TQ_{\Gamma_1}y^{(N)}-(x^{(N)})^TQ_{\Gamma_1}y_0\\
&+\breve{x}_0^T\Gamma_1^TQ\Gamma_1y_0+\breve{u}_iRu_i\big]dt+ \sum_{i=1}^N \mathbb{E}\big[\breve{x}_i^T(T)Hy_i(T)-(\breve{x}^{(N)}(T))^T{H_{\hat{\Gamma}}}y^{(N)}(T)\cr
&- \breve{x}_0^T(T)H_{\hat{\Gamma}_1}^Ty^{(N)}(T)-[\breve{x}^{(N)}(T)]^TH_{\hat{\Gamma}_1} y_0(T)+\breve{x}_0^T(T)\hat{\Gamma}_1^TH\hat{\Gamma}_1y_0(T)\big],
\label{eq8b} \nonumber
\end{align}
\begin{align}
I_2=&\sum_{i=1}^N \mathbb{E}\int_0^T\big[|y_i|_Q^2-|y^{(N)}|^2_{Q_{\Gamma}}-2\Gamma y_0^TQ_{\Gamma_1}^Ty^{(N)}+y_0^T\Gamma_1^TQ\Gamma_1y_0+|{u}_i|^2_R\big]dt\\
&+ \sum_{i=1}^N \mathbb{E}\big[|y_i(T)|_H^2-|y^{(N)}(T)|^2_{H_{\hat{\Gamma}}}
-2( y_0(T))^TH_{\hat{\Gamma}_1}^Ty^{(N)}(T)+|y_0(T)|^2_{\hat{\Gamma}_1^TH\hat{\Gamma}_1}\big].\nonumber
\end{align}
Let $\{\breve{p}_i, \breve{q}_i^j,i,j=0,1,\cdots,N\}$ be a set of solutions to (\ref{eq3}).
Then, by It\^{o}'s formula, we obtain 
\begin{align*}
	&\sum_{i=1}^N\mathbb{E}\big[\langle\hat{\Gamma}^T_1H (\hat{\Gamma}-I)\breve{x}^{(N)}(T)
	+\hat{\Gamma}_1^TH\hat{\Gamma}_1\breve{x}_0^T(T), z_0(T)\rangle\big]\cr
	= &\sum_{i=1}^N \mathbb{E}\int_0^T \Big\{\big\langle -[F\breve{p}^{(N)}+\bar{F}\breve{q}^{(N)}-\Gamma_1^TQ((I-\Gamma) \breve{x}^{(N)}-\Gamma_1\breve{x}_0)],z_0\big\rangle\cr
	&+\langle {G}_0^T\breve{p}_0+\bar{G}_0^T\breve{q}^0_0,z_i\big\rangle\Big\} dt,
\end{align*}
and
\begin{align*}
	&\sum_{i=1}^N\mathbb{E}[\langle H\breve{x}_i(T)-H_{\hat{\Gamma}}\breve{x}^{(N)}(T)+(\hat{\Gamma}-I)^TH\hat{\Gamma}_1 \breve{x}_0(T),z_i(T)\rangle]\cr
	= &\sum_{i=1}^N \mathbb{E}\int_0^T \Big\{\big\langle -\big[Q\breve{x}_i-Q_{\Gamma}\breve{x}^{(N)}
	+(\Gamma-I)^TQ\Gamma_1\breve{x}_0+{G}_0^T\breve{p}_0+\bar{G}_0^T\breve{q}^0_0\big],z_i\big\rangle\cr
	&+\langle F\breve{p}^{(N)}+\bar{F}\breve{q}^{(N)},z_0\rangle+\langle B^T\breve{p}_i+D^T\breve{q}_i^i,u_i\rangle\Big\} dt,
\end{align*}
where
the second equation holds since $\sum_{i=1}^N\mathbb{E}\langle G^T\breve{p}^{(N)},z_i\rangle=\sum_{i=1}^N\mathbb{E}\langle \breve{p}_i,Gz^{(N)}\rangle $
and  $\sum_{i=1}^N\mathbb{E}\langle \bar{G}^T\breve{q}^{(N)},z_i\rangle=\sum_{i=1}^N\mathbb{E}\langle \breve{q}_i^i,\bar{G}z^{(N)}\rangle .$
From the above equations and (\ref{eq8a}), 
\begin{align}\label{eqA3}
	I_1=&\frac{1}{N}\sum_{i=1}^N \mathbb{E}\int_0^T\big[\big\langle Q\breve{x}_i-{Q_{\Gamma}}\breve{x}^{(N)}+(\Gamma-I)^TQ\Gamma \breve{x}_0,z_i\big\rangle+\langle\Gamma_1^T Q(\Gamma-I)\breve{x}^{(N)} +\Gamma_1^TQ\Gamma_1\breve{x}_0,z_0\rangle\cr
	&+ \langle R\breve{u}_i,u_i\rangle\big]dt+ \frac{1}{N}\sum_{i=1}^N \mathbb{E}\big[\big\langle H\breve{x}_i(T)-{H_{\hat{\Gamma}}}\breve{x}^{(N)}(T)+(\hat{\Gamma}-I)^TH\hat{\Gamma} x_0(T),z_i(T)\big\rangle\cr
	&+\big\langle\hat{\Gamma}^T_1H (\hat{\Gamma}-I)\breve{x}^{(N)}(T)
	+\hat{\Gamma}_1^TH\hat{\Gamma}_1\breve{x}_0^T(T), z_0(T)\big\rangle\big]\cr
	=&\frac{1}{N}\sum_{i=1}^N \mathbb{E}\int_0^T[\langle R\breve{u}_i+B^T\breve{p}_i+D^T\breve{q}_i^i,u_i\rangle]dt.
\end{align}
 From (\ref{eqA1}), $\breve{u}$ is a minimizer to (P1) if and only if
$I_2\geq 0$ and $I_1=0$.  Indeed, if $I_2\geq 0$  does not hold, then the problem (P1) is ill-posed (see, e.g., \cite{SLY16}). By \cite[Proposition 3.1]{WZ21}, $I_2\geq 0$  if and only if (P1) is convex.
By (\ref{eqA3}),  $I_1=0 $ is equivalent to
$R\breve{u}_i+B^T\breve{p}_i+D^T\breve{q}_i^i=0, \ i=1,\cdots,N.$
Thus, we have the optimality system (\ref{eq59}), which implies that (\ref{eq59}) admits a solution $(\check{x}_i,\check{p}_i,\check{q}_{i}^{j}, i,j=1,\cdots,N)$.

On other hand, if the equation system (\ref{eq3}) admits a solution $\{x_i,p_i,q_i^j,i,j=1,\cdots,N\}$.
Let $\check{u}_i$ satisfy (\ref{eq3-b}).  If (P1) is convex, then by (\ref{eqA3}),
$\check{u}_i$
is an optimal control to Problem (P1).
\hfill{$\Box$}

To prove Theorem \ref{thm4.2}, we provide two lemmas.

\begin{lemma}\label{lem2a}
Assume that (A1)-(A4) hold. Then, the following holds:
\begin{align}
	&\sup_{0\leq t\leq T} \mathbb{E}\big[|\bar{x}^{(N)}-\bar{x}|^2 + |\bar{p}^{(N)}-\mathbb{E}_{\mathcal{F}^0}[\bar{p}_i]|^2+ |\bar{q}^{(N)}-\mathbb{E}_{\mathcal{F}^0}[\bar{q}^i_i]|^2\big]=O(\frac{1}{N}),
\end{align}
\end{lemma}
where 
$\bar{p}^{(N)}=\frac{1}{N}\sum_{i=1}^N\bar{p}_i$ and $\bar{q}^{(N)}=\frac{1}{N}\sum_{i=1}^N\bar{q}_i^i$.

\emph{Proof.} 
After applying $u_i^*$, $i=0,\cdots,N$, we have
\begin{align}\label{eq40a}
d\bar{x}_i=&\big[(A-B\Upsilon^{\dag}\Psi_2)\bar{x}_i+({G}-B\Upsilon^{\dag}\Psi_3)\bar{x}-B\Upsilon^{\dag}B^T\varphi+(F-B\Upsilon^{\dag}\Psi_1)\bar{x}_0\big]dt\\
&+\big[({C}-D\Upsilon^{\dag}\Psi_2)\bar{x}_i+(\bar{G}-D\Upsilon^{\dag}\Psi_3)\bar{x}-D\Upsilon^{\dag}B\varphi+(\bar{F}-D\Upsilon^{\dag}\Psi_1)\bar{x}_0\big]dW_i.\nonumber
\end{align}
By Assumption (A4), $\mathbb{E}\int_0^T|u_0^*|^2dt\leq c_1$. Then, it leads to $\mathbb{E}\int_0^T|x_0^*|^2dt\leq c_2$. By (\ref{eq7}), $\max_{0\leq t\leq T} \mathbb{E}[|\bar{x}(t)|^2]\leq c_3$. This further gives that $\sup_{0\leq t\leq T} \mathbb{E}[|\bar{x}_i(t)|^2]\leq c_4.$
By (\ref{eq40a}) and (\ref{eq7}), we obtain
\begin{align*}
d(\bar{x}^{(N)}&-\bar{x})=(A-B\Upsilon^{\dag}\Psi_2)
(\bar{x}^{(N)}-\bar{x})dt\cr
&+\frac{1}{N}\sum_{i=1}^N\big[({C}-D\Upsilon^{\dag}\Psi_2)\bar{x}_i+(\bar{G}-D\Upsilon^{\dag}\Psi_3)\bar{x}-D\Upsilon^{\dag}B\varphi+(\bar{F}-D\Upsilon^{\dag}\Psi_1)\bar{x}_0\big]dW_i,
\end{align*}
which gives
\begin{align*}
\bar{x}^{(N)}(t)&-\bar{x}(t)=\Xi(t,0)[\bar{x}^{(N)}(0)-\bar{x}(0)]\cr
&+\frac{1}{N}\sum_{i=1}^N\int_0^t\Xi(t,s)\big[({C}-D\Upsilon^{\dag}\Psi_2)\bar{x}_i+(\bar{G}-D\Upsilon^{\dag}\Psi_3)\bar{x}-D\Upsilon^{\dag}B\varphi+(\bar{F}-D\Upsilon^{\dag}\Psi_1)\bar{x}_0\big]dW_i(s).
\end{align*}
Here, $\Xi(t,s)$ satisfies $\frac{d\Xi(t,s)}{dt}=(A-B\Upsilon^{\dag}\Psi_2)\Xi(t,s),\ \Xi(s,s)=I.$
By (A1), we further have
\begin{align}\label{eq41a}
&\mathbb{E}|\bar{x}^{(N)}(t)-\bar{x}(t)|^2\\
\leq&\big|\Xi(t,0)\big|^2\mathbb{E}|\bar{x}^{(N)}(0)-\bar{x}(0)|^2
+\frac{1}{N^2}\sum_{i=1}^N\int_0^t c_1\big|\Xi(t,s)\big|^2\max_{1\leq i\leq N}\mathbb{E}\big(|\bar{x}_i|^2+|\bar{x}|^2+|\varphi|^2+|\bar{x}_0|^2)ds\cr
\leq& \frac{1}{N}\Big\{\big|\Xi(t,0)\big|^2 \max_{1\leq i\leq N} \big[\mathbb{E}|x_{i0}|^2+ c_2\sup_{0\leq t\leq T}\mathbb{E}\big(|\bar{x}_i|^2+|\bar{x}|^2+|\varphi|^2+|\bar{x}_0|^2)\big]\Big\}=O(\frac{1}{N}).\nonumber
\end{align}
Note that $\bar{p}_i=P\bar{x}_i+\bar{P}\bar{x}+{P}_0\bar{x}_0+\varphi$. Then, we have
$$
\begin{aligned}
\sup_{0\leq t\leq T} \mathbb{E}\big[ |\bar{p}^{(N)}(t)-\mathbb{E}_{\mathcal{F}^0}[\bar{p}_i(t)]|^2\big]=\sup_{0\leq t\leq T} \mathbb{E}\big[ |P(\bar{x}^{(N)}(t)-\bar{x}(t))|^2\big]=O({1}/{N}).
\end{aligned}
$$
From (\ref{eq14a}), (\ref{eq10e}) and (\ref{eq41a}), we obtain
$$\begin{aligned}
\sup_{0\leq t\leq T} \mathbb{E}\big[ |\bar{q}^{(N)}(t)-\mathbb{E}_{\mathcal{F}^0}[\bar{q}^i_i(t)]|^2=\sup_{0\leq t\leq T} \mathbb{E}\big[ |P({C}-D\Upsilon^{\dag}\Psi_2)(\bar{x}^{(N)}(t)-\bar{x}(t))|^2\big]=O({1}/{N}).
\end{aligned}$$
\rightline{$\Box$}

\begin{lemma}\label{lem2}
Assume that (A1)-(A4) hold. Then, the following holds:
\begin{equation*}
	\begin{aligned}
		&\sup_{0\leq t\leq T} \mathbb{E}\big[|{x}_0^*(t)-\bar{x}_0(t)|^2+|{x}_*^{(N)}(t)-\bar{x}(t)| ^2\big]=O(\frac{1}{N}),\cr
		&\sup_{0\leq t\leq T}   \mathbb{E}|{x}_i^*(t)-\bar{x}_i(t)|^2=O(\frac{1}{N}),
	\end{aligned}
\end{equation*}	
where $x_i^*,i=1,\cdots,N$ is the realized state under the control $u_i^*,i=1,\cdots,N$.
\end{lemma}
\emph{Proof.} By (\ref{eq7}) and (\ref{eq23-a}), it can be verified that
$\sup_{0\leq t\leq T} \mathbb{E}\big[|\bar{x}^{(N)}-\bar{x}|^2$ and
$\max_{1\leq i\leq N}\mathbb{E}\int_0^T(|x_i^*|^2+|u_i^*|^2)dt\leq c_3.$
From (\ref{eq7}), 
we have
$$
\begin{aligned}
d({x}_0^*-\bar{x}_0)=& [A_0({x}_0^*-\bar{x}_0)+G_0({x}^{(N)}_*-\bar{x})]dt+[C_0({x}_0^*-\bar{x}_0)+\bar{G}_0({x}^{(N)}_*-\bar{x})]dW_0,\ {x}_0^*(0)-\bar{x}_0(0)=0,\\
d({x}^{(N)}_*-\bar{x})=&({A}+G)
({x}_*^{(N)}-\bar{x})dt+F({x}_0^*-\bar{x}_0)-B\Upsilon^{\dag}\Psi_2(\bar{x}^{(N)}-\bar{x})\cr
&+\frac{1}{N}\sum_{j=1}^N(C{x}_j^*+Du_j^*+\bar{G}{x}^{(N)}_*+\bar{F}{x}_0^*)dW_j\  {x}_*^{(N)}(0)-\bar{x}(0)=\frac{1}{N}\sum_{i=1}^N\xi_i-\bar{\xi}.
\end{aligned}
$$
Similar to (\ref{eq41a}), we  have
\begin{align}\label{eq44b}
\sup_{0\leq t\leq T} \mathbb{E}\big[|{x}_0^*(t)-\bar{x}_0(t)|^2+|{x}_*^{(N)}(t)-\bar{x}(t)| ^2\big]=O(\frac{1}{N})
\end{align}
From (\ref{eq23-a}) and (\ref{eq40a}),
\begin{align*}
d({x}_i^*-\bar{x}_i)=[A({x}_i^*-\bar{x}_i)+G({x}_*^{(N)}-\bar{x})+F({x}_0^*-\bar{x}_0)]dt
+[C({x}_i^*-\bar{x}_i)+\bar{G}({x}_*^{(N)}-\bar{x})+\bar{F}({x}_0^*-\bar{x}_0)]dW_i,
\end{align*}
with ${x}_i^*(0)-\bar{x}_i(0)=0$.
Let $\Xi_i(t)$ be the solution to the following stochastic differential euquation:
$d\Xi_i(t)=A\Xi_i(t)dt+C\Xi_i(t)dW_i(t),\ \Xi_i(0)=I.$
Then, one can obtain
\begin{align*}
{x}_i^*(t)-\bar{x}_i(t)=&\int_0^t\Xi_i(t)\Xi_i^{-1}(s)\big[G({x}_*^{(N)}(s)-\bar{x}(s))+F({x}_0^*(s)-\bar{x}_0(s))\big]ds\cr
&+\int_0^t\Xi_i(t)\Xi_i^{-1}(s)\big[\bar{G}({x}_*^{(N)}(s)-\bar{x}(s))+\bar{F}({x}_0^*(s)-\bar{x}_0(s))\big]dW_i(s).
\end{align*}
Note that $\mathbb{E}\int_0^T|\Xi_i^T(t)\Xi_i(t)|dt<c$. From (\ref{eq44b}), we have
\begin{align*}
\mathbb{E}| {x}_i^*(t)-\bar{x}_i(t)|^2\leq &2T\mathbb{E}\int_0^t\big|\Xi_i(t)\Xi_i^{-1}(s)\big|^2\big|G({x}_*^{(N)}(s)-\bar{x}(s))+F({x}_0^*(s)-\bar{x}_0(s))\big|^2ds\cr	
&+2\mathbb{E}\int_0^t\big|\Xi_i(t)\Xi_i^{-1}(s)\big|^2\big|\bar{G}({x}_*^{(N)}(s)-\bar{x}(s))+\bar{F}({x}_0^*(s)-\bar{x}_0(s))\big|^2ds=O(\frac{1}{N}).
\end{align*}
This with (\ref{eq44b}) completes the proof. \hfill{$\Box$}

\emph{Proof of Theorem \ref{thm4.2}.}
\emph{(For followers).} 
We first prove that for $u\in \mathcal{U}_c$,  $J_{\rm soc}(u)< \infty$ implies that
$\mathbb{E}\int_0^T(|x_i|^2+|u_i|^2)dt<\infty,$ for all $i=1,\cdots,N$.
In views of  (A2), by \cite{SLY16} we have
$$\delta_0 \sum_{i=1}^N\mathbb{E}\int_0^T |u_i|^2dt-c_0\leq J_{\rm soc}(u)< \infty,$$
which implies
$\sum_{i=1}^N\mathbb{E}\int_0^T|u_i|^2dt<c_1.$
By (\ref{eq1}) and Schwarz's inequality \cite{YZ99},
\begin{align*}
\mathbb{E}|x_i(t)|^2&\leq c_2\mathbb{E}\int_0^t|x^{(N)}(\tau)|^2d\tau+c_3\leq\frac{c_2}{N}\mathbb{E}\int_0^t\sum_{j=1}^N|x_j(\tau)|^2d\tau+c_3.
\end{align*}
By Gronwall's inequality,
we have $\sum_{j=1}^N\mathbb{E}|x_j(t)|^2\leq Nc_3e^{c_2t}\leq Nc_3e^{c_2T}.$

Let $\tilde{u}_i=u_i-{u}^*_i$, $\tilde{x}_i=x_i-{x}_i^*$, $\tilde{x}^{(N)}=\frac{1}{N}\sum_{i=1}^N \tilde{x}_i$ and $\tilde{x}_0=x_0-{x}_0^*$. Then, by (\ref{eq1}) and (\ref{eq23-a}), we get
\begin{equation}\label{eq32}
\left\{\begin{aligned}
d\tilde{x}_0 =&(A_0\tilde{x}_0 +{G}_0\tilde{x}^{(N)} )dt++(C_0\tilde{x}_0 +\bar{G}_0\tilde{x}^{(N)} )dW_0,\ \tilde{x}_0(0)=0,\\
d\tilde{x}_i =&(A\tilde{x}_i +{G}\tilde{x}^{(N)}+F\tilde{x}_0 +B\tilde{u}_i )dt+(C\tilde{x}_i +\bar{G}\tilde{x}^{(N)} +\tilde{F}\tilde{x}_0+D\tilde{u}_i )dW_i,\ \tilde{x}_i(0)=0.
\end{aligned}
\right.
\end{equation}
From (\ref{eq3}), we have
$J^{(N)}_{\rm soc}(u_0^*,u)=\frac{1}{N}\sum_{i=1}^N({J}_i(u_0^*,{u}^*)+\tilde{J}_i(u_0^*,\tilde{u})+\mathcal{I}_i),$
where
\begin{align*}
&\tilde{J}_i(u_0^*,\tilde{u})\stackrel{\Delta}{=}\mathbb{E}\int_0^T\big[|\tilde{x}_i-\Gamma \tilde{x}^{(N)}-\Gamma _1 \tilde{x}_0|^2_Q+|\tilde{u}_i|^2_{R}\big]dt\cr
&\qquad\qquad\quad +\mathbb{E}|\tilde{x}_i(T)-\hat{\Gamma} \tilde{x}^{(N)}(T)-\hat{\Gamma} _1 \tilde{x}_0(T)|^2_H,\\
&\mathcal{I}_i=2\mathbb{E}\int_0^T\big[\big({x}_i^*-\Gamma {x}_*^{(N)}-\Gamma _1 {x}^*_0\big)^TQ\big(\tilde{x}_i
-\Gamma \tilde{x}^{(N)}-\Gamma _1 \tilde{x}_0\big)+\tilde{u}^T_iR{u}_i^*\big]dt\cr
&\qquad+\mathbb{E}\big[\big({x}_i^*(T)-\hat{\Gamma} {x}_*^{(N)}(T)-\hat{\Gamma} _1 {x}^*_0(T)\big)^TH\big(\tilde{x}_i(T)-\hat{\Gamma} \tilde{x}^{(N)}(T)-\hat{\Gamma} _1 \tilde{x}_0(T)\big)\big].
\end{align*}
By (\ref{eq32}) and It\^{o}'s formula,
we obtain
$$  \begin{aligned}
&N\mathbb{E}\big[  \tilde{x}_0^T(T)\big(-H_{\bar{\Gamma}_1} ^T\bar{x}(T)+\bar{\Gamma}_1^TH\bar{\Gamma}_1\bar{x}_0(T)\big)\big]
\cr
=&\sum_{i=1}^N\mathbb{E}\int_{0}^{T}\Big\{\tilde{x}_i^T (G_0^T\bar{p}_0+\bar{G}_0^T\bar{q}^0_0)- \tilde{x}_0^T\big[ F^T\bar{p}^{(N)}+\bar{F}^T\bar{q}^{(N)}+\Gamma^T_1Q((\Gamma-I) \bar{x}+\Gamma_1\bar{x}_0)
&+F^T\Pi(\bar{x}-\bar{x}^{(N)})+\bar{F}^T(\mathbb{E}_{\mathcal{F}_0}[\bar{q}_i^i]-\bar{q}^{(N)})\big]\Big\}dt,
\end{aligned}$$
$$  \begin{aligned}
&\sum_{i=1}^N\mathbb{E}\big[\tilde{x}_i^T(T)\big(H\bar{x}_i(T)-H_{\hat{\Gamma}}\bar{x}(T)-H_{\hat{\Gamma}_1} {x}_0^*(T)\big)\big]=\sum_{i=1}^N\mathbb{E}\big[  \tilde{x}_i^T(T)\bar{p}_i(T))\big]\cr
=&\mathbb{E}\int_{0}^{T}\sum_{i=1}^N\Big\{- \tilde{x}_i^T\big[G_0^T\bar{p}_0+\bar{G}_0^T\bar{q}_0^0+ Q\bar{x}_i-Q_{\Gamma} \bar{x}+(\Gamma-I)^TQ\Gamma_1x_0^*\big]+\tilde{x}_0^T(F^T\bar{p}_i+\bar{F}^T\bar{q}_i^i)-\tilde{u}_i^TR{u}_i^*\Big\}dt\cr
&+\sum_{i=1}^N\mathbb{E}\int_{0}^{T}\tilde{x}_i^T\big[G^T(\bar{p}^{(N)}
-\mathbb{E}_{\mathcal{F}^0}[\bar{p}_i])+\bar{G}^T(\bar{q}^{(N)}-\mathbb{E}_{\mathcal{F}^0}[\bar{q}^i_i])\big]dt.
\end{aligned}$$
From this and direct computations, one can obtain 
\begin{equation*}
\begin{aligned}
	\frac{1}{N}\sum_{i=1}^N \mathcal{I}_i=&\frac{1}{N} \sum_{i=1}^N2\mathbb{E}\Big\{\int_{0}^{T}\tilde{x}_i^T\big[Q(x_i^*-\bar{x}_i))+Q_{\Gamma}({x}_*^{(N)}-\bar{x})+G^T(\bar{p}^{(N)}-\mathbb{E}_{\mathcal{F}^0}\big[\bar{p}_i])+F^T\Pi(\bar{x}-\bar{x}^{(N)})\cr
	&+(\bar{G}-\bar{F})^T(\bar{q}^{(N)}-\mathbb{E}_{\mathcal{F}^0}[\bar{q}^i_i])\big]dt+\big[\tilde{x}_i^T(T)\big(H({x}^*_i(T)-\bar{x}_i(T))-H_{\hat{\Gamma}}({x}^{(N)}_*(T)-\bar{x}(T))\big]\Big\}\cr
	\leq &\frac{c}{N} \sum_{i=1}^N\Big[\mathbb{E}\int_{0}^{T}|\tilde{x}_i|^2dt\Big]^{1/2}\cdot  \Big[\mathbb{E}\int_{0}^{T}\big(|x_i^*-\bar{x}_i|^2+|{x}_*^{(N)}-\bar{x}|^2+|\bar{x}^{(N)}-\bar{x}|^2\cr
	&+|\bar{p}^{(N)}-\mathbb{E}_{\mathcal{F}^0}\big[\bar{p}_i]|^2+|\bar{q}^{(N)}-\mathbb{E}_{\mathcal{F}^0}[\bar{q}^i_i|^2]\big)dt\Big]^{1/2}+ O\big(\frac{1}{\sqrt{N}}\big)
	\leq O\big({1}/{\sqrt{N}}\big)=\epsilon_1.
\end{aligned}
\end{equation*}
Note that by (A2), $\sum_{i=1}^N\tilde{J}_i(\tilde{u},u_0^*)\geq0$.
Then, we have
${J}_{\rm soc}({u}^*,u_0^*)\leq J_{\rm soc}(u,u_0^*)+\epsilon_1.$

(\emph{For the leader}).
By (\ref{eq8-a}) and Schwarz's inequality,  we have
\begin{align}\label{eq36-a}
J_0({u}_0^*,{u}^*)=&\mathbb{E}\int_{0}^{T}\big[|\bar{x}_0^*- \Gamma_0\bar{x}+{x}_0^*-\bar{x}_0^*+\Gamma_0({x}^{(N)}_*-\bar{x})|_{Q_0}^{2}+|{u}_0^*|_{R_0}^{2}\big]dt\\
&+\mathbb{E}\big[|\bar{x}_0^*(T)- \hat{\Gamma}_0\bar{x}(T)+{x}_0^*(T)-\bar{x}_0^*(T)+\hat{\Gamma}_0({x}^{(N)}_*(T)-\bar{x}(T))|_{H_0}^{2}\big]dt\cr
\leq& \bar{J}_0({u}_0^*,{u}^*)+2\mathbb{E}[|{x}_0^*-\bar{x}_0^*|_{Q_0}^{2}+|\Gamma_0({x}_*^{(N)}-\bar{x})|_{Q_0}^{2}]\big]dt\cr
&+2\mathbb{E}\big[|{x}_0^*(T)-\bar{x}_0^*(T)|_{Q_0}^{2}+|\hat{\Gamma}_0({x}_*^{(N)}(T)-\bar{x}(T))|_{H_0}^{2}\big]\cr
&+C\sup_{0\leq t\leq T}\big(\mathbb{E}[|{x}_0^*(t)-\bar{x}_0^*(t)|^{2}+|{x}_*^{(N)}(t)-\bar{x}(t)|^2]\big)^{1/2}
\cr
\leq &\bar{J}_0({u}_0^*,{u}^*)+O(1/\sqrt{N}).\nonumber
\end{align}
It follows from Theorem \ref{thm4.0} that
$  \bar{J}_0({u}_0^*,{u}^*) \leq \bar{J}_0(u_0,{u}^*).
$
This together with (\ref{eq36-a}) implies
\begin{equation}\label{eq41c}
J_0({u}^*_0,{u}^*({u}^*_0))\leq  \bar{J}_0(u_0,u^*({u}_0))+O(1/\sqrt{N}),
\end{equation}
for any ${u}_0\in {\cal U}_0$. From (\ref{eq8-a}), we obtain
\begin{align*}
\bar{J}_0({u}_0,{u}^*)=&\mathbb{E}\int_{0}^{T}\big[|x_0- \Gamma_0{x}^{(N)}_*+\bar{x}_0^*-x_0+\Gamma_0({x}_*^{(N)}-\bar{x})|_{Q_0}^{2}+|u_0|_{R_0}^{2}\big]dt\cr
&+\mathbb{E}\big[|{x}_0(T)- \bar{\Gamma}_0{x}^{(N)}_*(T)+\bar{x}_0(T)-{x}_0(T)+\bar{\Gamma}_0({x}^{(N)}_*(T)-\bar{x}(T))|_{H_0}^{2}\big]dt\cr
\leq &{J}_0({u}_0,{u}^*)+O(1/\sqrt{N}),
\end{align*}
which  with (\ref{eq41c}) gives
$J_0({u}^*_0,{u}^*({u}^*_0))\leq {J}_0({u}_0,{u}^*({u}_0))+\varepsilon_2,$
where $\varepsilon_2=O(1/\sqrt{N})$.
\hfill{$\Box$}

\section{Proof of Theorem \ref{thm3.5}}\label{appb}
\def\theequation{B.\arabic{equation}}
\setcounter{equation}{0}
To prove Theorem \ref{thm3.5}, we first give a lemma.
Consider an MF-type problem: optimize the cost functional
\begin{align}\label{eq39}
\mathcal{J}_i(u_i)=&\mathbb{E}\int_0^T\big(|\bar{x}_i-{\Gamma}\mathbb{E}_{\mathcal{F}^0}[\bar{x}_i]-\Gamma_1\bar{x}_0|^2_Q
+|u_i|^2_R\big)dt+\mathbb{E}\big[|\bar{x}_i(T)-\hat{\Gamma}\mathbb{E}_{\mathcal{F}^0}[\bar{x}_i(T)]-\hat{\Gamma}_1\bar{x}_0(T)|^2_H\big]
\end{align}
subject to ($\bar{x}_i(0)=\xi_i$)
\begin{equation}
  \label{eq40-a}
  \left\{
\begin{aligned}
d\bar{x}_0=&(A_0\bar{x}_0+B_0{u}_0^*+G_0\mathbb{E}_{\mathcal{F}^0}[\bar{x}_i])dt+(C_0\bar{x}_0+D_0{u}_0^*+\bar{G}_0\mathbb{E}_{\mathcal{F}^0}[\bar{x}_i])dW_0,\cr
d\bar{x}_i=&(A\bar{x}_i+B{u}_i+G\mathbb{E}_{\mathcal{F}^0}[\bar{x}_i]+F\bar{x}_0)dt+(C\bar{x}_i+D{u}_i+\bar{G}\mathbb{E}_{\mathcal{F}^0}[\bar{x}_i]+\bar{F}\bar{x}_0)dW_i.
\end{aligned}
\right.
\end{equation}

\begin{lemma}\label{lem1}
Assume (A1) and (A4) hold.
For Problem (\ref{eq39})-(\ref{eq40-a}), the optimal control $u_i^*$ is given by (\ref{eq10e}),
and the corresponding optimal cost is $\mathbb{E}[|\xi_i|^2_{P(0)}+|\bar{\xi}_0|^2_{K(0)}+2\varphi^T(0)\bar{x}_0]+s_T$. 

\end{lemma}

\emph{Proof.} Note that $\mathbb{E}_{\mathcal{F}^0}[\bar{x}_i]=\bar{x}$ satisfies
$d\bar{x}= \big[(A+G)\bar{x}+B\bar{u}+F\bar{x}_0\big]dt,$
where $\bar{u}=\mathbb{E}_{\mathcal{F}^0}[{u}_i]$. 
 Applying It\^{o}'s formula to $\|\bar{x}_i\|^2_P$ yields
\begin{equation}\label{eqb4}
\begin{aligned}
 & \mathbb{E}|\bar{x}_i(T)|^2_H-\mathbb{E}|\bar{x}_i(0)|^2_{P(0)}=\mathbb{E}\int_0^T\Big[\bar{x}_i^T(\dot{P}+PA+A^TP)\bar{x}_i+\bar{x}^T(PG+G^TP)\bar{x}\cr
  &+2\bar{x}_i^TPBu_i+2\bar{x}_i^TPF\bar{x}_0+|C\bar{x}_i+D{u}_i+\bar{G}\bar{x}+\bar{F}\bar{x}_0|^2_{P}\Big]dt.
\end{aligned}
\end{equation}
Also, applying It\^{o}'s formula, we have
\begin{equation}
\begin{aligned}
 & -\mathbb{E}|\bar{x}(T)|^2_{H_{\hat{\Gamma}}}-\mathbb{E}|\bar{x}(0)|^2_{\bar{P}(0)}=\mathbb{E}\int_0^T\Big[\bar{x}^T(\dot{\bar{P}}+\bar{P}(A+G)+(A+G)^T\bar{P})\bar{x}+2\bar{x}^T\bar{P}(B\bar{u}+F\bar{x}_0)\Big]dt,
\end{aligned}
\end{equation}
\begin{equation}
\begin{aligned}
 & \mathbb{E}|\bar{x}_0(T)|^2_{\hat{\Gamma}_1^TH\hat{\Gamma}_1}-\mathbb{E}|\bar{x}_0(0)|^2_{K(0)}=\mathbb{E}\int_0^T\Big[\bar{x}_0^T(\dot{K}+KA_0+A_0^TK)\bar{x}_0\cr
  &+2\bar{x}_0^TK(B_0{u}_0^*+G_0\bar{x})+|C_0\bar{x}_0+D_0{u}_0^*+\bar{G}_0\bar{x}|^2_{K}\Big]dt,
\end{aligned}
\end{equation}
\begin{equation}
\begin{aligned}
 & \mathbb{E}\big[-\bar{x}_0(T)H^T_{\hat{\Gamma}_1}\bar{x}(T)-\bar{x}_0^T(0)P_0(0)\bar{x}(0)\big]=\mathbb{E}\int_0^T\Big[\bar{x}_0^T\big(\dot{P}_0+P_0(A+G)+A^TP_0\big)\bar{x}\cr
  &+(G_0\bar{x}+B_0u_0^*)^TP_0\bar{x}+\bar{x}_0^TP_0(B\bar{u}+{F}\bar{x}_0)\Big]dt,
\end{aligned}
\end{equation}
\begin{equation}
\begin{aligned}
 & \mathbb{E}\big[\bar{x}^T(T)\varphi(T)-\bar{x}^T(0)\varphi(0)\big]=\mathbb{E}\int_0^T\Big[(B\bar{u}+F\bar{x}_0)^T\varphi\cr
  &-\bar{x}^T\big(-(\Psi_2+\Psi_3)^T\Upsilon^{\dag}B^T\varphi+G_0^T\varphi_0+\bar{G}_0^T\zeta_0+P_0B_0u_0^*\big)\Big]dt,
\end{aligned}
\end{equation}
and
\begin{equation}\label{eqb9}
\begin{aligned}
 & \mathbb{E}\big[\bar{x}^T_0(T)\varphi_0(T)-\bar{x}^T_0(0)\varphi_0(0)\big]=\mathbb{E}\int_0^T\Big[(G_0\bar{x}+B_0{u}_0)^T\varphi_0+(\bar{G}_0\bar{x}+D_0{u}_0^*)^T\zeta_0\cr
  &-\bar{x}^T_0\big((F-B\Upsilon^{\dag}\Psi_1)^T\varphi+(C_0^TKD_0+KB_0)u_0^*\big)\Big]dt.
\end{aligned}
\end{equation}
By (\ref{eqb4})-(\ref{eqb9}), we obtain
\begin{align*}
\mathcal{J}_i(u_i)
=&\mathbb{E}\int_0^T\big(|\bar{x}_i|_Q^2-|\bar{x}|^2_{Q_{\Gamma}}-2\bar{x}Q_{\Gamma_1}\bar{x}_0+|\bar{x}_0|^2_{\Gamma_1^TQ\Gamma_1}
+|u_i|^2_R\big)dt\cr
&
+\mathbb{E}\big[|\bar{x}_i(T)|_H^2-|\bar{x}(T)|^2_{H_{\hat{\Gamma}}}-2\bar{x}(T)H_{\hat{\Gamma}_1}\bar{x}_0(T)+|\bar{x}_0(T)|^2_{\hat{\Gamma}_1^TH\hat{\Gamma}_1}\big]\cr
=&\mathbb{E}[|x_{i0}|^2_{P(0)}+|\bar{x}(0)|^2_{\bar{P}(0)}+|\bar{x}_0(0)|^2_{K(0)}+2\bar{x}_0^T(0)P_0(0)\bar{x}(0)+2\varphi^T(0)\bar{x}(0)+2\varphi_0^T(0)\bar{x}_0(0)]\cr
&+\mathbb{E}\int_0^T\Big[ \bar{x}_i^T\Psi_2^T\Upsilon^{\dag}\Psi_2\bar{x}_i+ \bar{x}^T(\Psi_2^T\Upsilon^{\dag}\Psi_3+\Psi_3^T\Upsilon^{\dag}\Psi_2+\Psi_3^T\Upsilon^{\dag}\Psi_3)\bar{x}
+\bar{x}_0^T\Psi_1^T\Upsilon^{\dag}\Psi_1\bar{x}_0\cr
&+\bar{x}_0^T(\Psi_2+\Psi_3)^T\Upsilon^{\dag}\Psi_1\bar{x}+2u_i^T(\Psi_2\bar{x}_i+\Psi_3\bar{x}+\Psi_1\bar{x}_0+B\varphi)+2\bar{x}^T(\Psi_2+\Psi_3)^T\Upsilon^{\dag}B^T\varphi\cr
&+2\bar{x}^T_0\Psi_1^T\Upsilon^{\dag}B^T\varphi+|u_i|^2_{\Upsilon}+(\bar{x}^T\bar{G}_0^TKD_0+\zeta_0^TD_0+\varphi_0^TB_0)u_0^*+u_0^*D_0KD_0u_0^*\Big]dt\cr
=&\mathbb{E}[|\xi_i|^2_{P(0)}+|\bar{\xi}|^2_{\bar{P}(0)}+|\bar{\xi}_0|^2_{K(0)}+2\bar{\xi}_0^TP_0(0)\bar{\xi}+2\varphi^T(0)\bar{\xi}+2\varphi^T_0(0)\bar{\xi}_0]\cr
&+\mathbb{E}\int_0^T\big[|u_i+\Upsilon^{\dag}(\Psi_2\bar{x}_i+\Psi_3\bar{x}+B^T\varphi+\Psi_1\bar{x}_0)|^2_{\Upsilon}\big]dt+s_T\cr
\geq & \mathbb{E}[|\xi_i|^2_{P(0)}+|\bar{\xi}|^2_{\bar{P}(0)}+|\bar{\xi}_0|^2_{K(0)}+2\bar{\xi}_0^TP_0(0)\bar{\xi}+2\varphi^T(0)\bar{\xi}+2\varphi^T_0(0)\bar{\xi}_0]+s_T.
\end{align*}
\rightline{$\Box$}
\emph{Proof of Theorem \ref{thm3.5}.}
Applying the control 
(\ref{eq27a}) into the social cost, it follows that 
\noindent
\begin{align*}
&J_{\rm soc}^{(N)}
(u^*,u_0^*)
=\frac{1}{N}\sum_{i=1}^N\mathbb{E}\Big\{\int_0^T\big[|\bar{x}_i-\Gamma \bar{x}-\Gamma_1\bar{x}_0+x_i^*-\bar{x}_i-\Gamma (x_*^{(N)}-\bar{x})-\Gamma_1(x_0^*-\bar{x}_0)|^2_Q\cr
&+ |\Upsilon^{\dag}(B^TP+D^TPC) \bar{x}_i+(B^TK+D^TP\bar{G})\bar{x}+B^T\varphi+D^TP\bar{F}\bar{x}_0)|^2_{R}\big]dt\cr
&+|\bar{x}_i(T)-\hat{\Gamma} \bar{x}(T)-\hat{\Gamma}_1\bar{x}_0(T)
+x_i^*(T)-\bar{x}_i(T)-\hat{\Gamma} (x_*^{(N)}(T)-\bar{x}(T))-\Gamma_1(x_0^*(T)-\bar{x}_0(T))|^2_H\Big\}.
\end{align*}
By Lemma \ref{lem2} and Schwarz's inequality, one can obtain
\begin{align*}
&\big|J_{\rm soc}^{(N)}(u^*,u_0^*)-\sum_{i=1}^N \mathcal{J}_i(u_i^*)  \big|\cr
\leq& \frac{1}{N}\sum_{i=1}^N\mathbb{E}\int_0^T\big[|x_i^*-\bar{x}_i|^2_Q+|\Gamma (x_*^{(N)}-\bar{x})|^2_Q+|\Gamma_1(x_0^*-\bar{x}_0)|^2_Q\big]dt+\frac{c}{N}\sum_{i=1}^N\sup_{0\leq t\leq T}\big(\mathbb{E}|x_i^*-\bar{x}_i|^2_Q\big)^{1/2}\cr
&+\frac{C}{N}\sum_{i=1}^N
\sup_{0\leq t\leq T}\big(\mathbb{E}|\Gamma (x_*^{(N)}-\bar{x})|^2_Q\big)^{1/2}+\frac{c}{N}\sum_{i=1}^N\sup_{0\leq t\leq T}\big(\mathbb{E}|\Gamma_1(x_0^*-\bar{x}_0)|^2_Q\big)^{1/2}
\leq O(\frac{1}{\sqrt{N}}).
\end{align*}
This together with Lemma \ref{lem1} leads to (\ref{eq426}).

(For the leader) By a similar argument with the proof of Theorem \ref{thm4.0}, one can obtain
$$\begin{aligned}
\bar{J}_{0}(u_0^*,u^*)=\mathbb{E}\Big\{{\xi}^T_0y_0(0)+\bar{\xi}^T\bar{y}(0)+\int_0^T\big[\big\langle R_0u_0^*+B_0^Ty_0+\bar{B}_1^T\bar{y}, u_0^*\big\rangle \big]dt\Big\}.
\end{aligned}$$
By (\ref{4.12}), we have $\lim_{N\to\infty}J_{0}(u_0^*,u^*)=\mathbb{E}\big[{\xi}^T_0y_0(0)+\bar{\xi}^T\bar{y}(0)\big] $. Thus, the theorem follows.
\hfill{$\Box$}

\section{Proof of Theorem  \ref{thm4.4}}
\def\theequation{C.\arabic{equation}}
\setcounter{equation}{0}

\emph{Proof of Theorem \ref{thm5.1}.}
{Suppose that $\{\breve{u}_i,i=1,\cdots,N\}$ is an optimal control of Problem (P3).}
Denote by $\breve{x}_i$ the state of player $i$ under the optimal control $\breve{u}_i$. For any $u_i\in L^2_{{\mathcal F}}(0, T; \mathbb{R}^r) $ and $\lambda\in \mathbb{R}\ (\lambda \not= 0)$, let $u_i^{\lambda}=\breve{u}_i+\lambda u_i$, $i=1,\cdots,N$. Denote by $x_0^{\lambda},x_i^{\lambda}$ the solution to the following perturbed  equation:
$$\left\{ \begin{aligned}
dx^{\lambda}_0=&\big[A_0x^{\lambda}_0+B_0(P_0x^{\lambda}_0+\bar{P}x_{\lambda}^{(N)})+G_0x_{\lambda}^{(N)})\big]dt+\big[C_0x^{\lambda}_0+D_0(P_0x_0^{\lambda}+\bar{P}x_{\lambda}^{(N)})+\bar{G}_0x_{\lambda}^{(N)}\big]dW_0,\cr
dx_i^{\lambda}=&\big(Ax_i^{\lambda}+B(\check{u}_i+\lambda u_i)+Gx_{\lambda}^{(N)}+Fx_0^{\lambda}\big)dt+\big(Cx_i^{\lambda}+Du_i^{\lambda}+\bar{G}x_{\lambda}^{(N)}+\bar{F}x_0^{\lambda}\big) dW_i,\cr
x_0^{\lambda}(0)&=\xi_{0},\ x_i^{\lambda}(0)=\xi_{i},\ i=1,2,\cdots,N,
\end{aligned}\right.$$
with $x_{\lambda}^{(N)}=\frac{1}{N}\sum_{i=1}^Nx^{\lambda}_i $. Let $z_i=(x_i^{\lambda}-\check{x}_i)/\lambda$. 
It can be verified that
$z_i$ satisfies 
\begin{equation*}
\left\{ \begin{aligned}
	dz_0=&\big[(A_0+B_0P_0)z_0+(G_0+B_0\bar{P})z^{(N)}\big]dt+\big[(C_0+D_0P_0)z_0+(\bar{G}_0+D_0\bar{P})z^{(N)}\big]dW_0,\  z_0(0)=0,\cr
	dz_i=&[A z_i+B u_i+Gz^{(N)}+Fz_0]dt+[C z_i+D u_i+\bar{G}z^{(N)}+\bar{F}z_0]dW_i,\  z_i(0)=0, \
\end{aligned}\right.
\end{equation*}
where $ i=1,2,\cdots,N$, and $z^{(N)}=\frac{1}{N}\sum_{i=1}^Nz_i$.
From (\ref{eq7-a}), we have
\begin{equation}\label{eq5-3}
J_{\rm  soc}^{(N)}(\breve{u}+\lambda u)-J^{(N)}_{\rm soc}(\breve{u})=2\lambda I_1+\lambda^2 I_2,
\end{equation}
where
\begin{align}\label{eq60}
I_1=&\frac{1}{N}\sum_{i=1}^N \mathbb{E}\int_0^T\big[\breve{x}_i^TQz_i-(\breve{x}^{(N)})^T{Q_{\Gamma}}z^{(N)}-\breve{x}_0^TQ_{\Gamma_1}{\hat{\Gamma}_1}^Tz^{(N)}-(\breve{x}^{(N)})^TQ_{\Gamma_1} z_0\\
&+\breve{x}_0^T\Gamma_1^TQ\Gamma_1z_0+\breve{u}_iRu_i\big]dt+ \sum_{i=1}^N \mathbb{E}\big[\breve{x}_i^T(T)Hz_i(T)-(\breve{x}^{(N)}(T))^T{H_{\hat{\Gamma}}}z^{(N)}(T)\cr
&- \breve{x}_0^T(T)H_{\hat{\Gamma}_1}^Tz^{(N)}(T)-[\breve{x}^{(N)}(T)]^TH_{\hat{\Gamma}_1} z_0(T)+\breve{x}_0^T(T)\hat{\Gamma}_1^TH\hat{\Gamma}_1z_0(T)\big],
\label{eq61} \nonumber
\end{align}
\begin{align}
I_2=&\frac{1}{N}\sum_{i=1}^N \mathbb{E}\int_0^T\big[|z_i|_Q^2-|z^{(N)}|^2_{Q_{\Gamma}}-2\Gamma z_0^TQ_{\Gamma_1}^Tz^{(N)}+z_0^T\Gamma_1^TQ\Gamma_1z_0+|{u}_i|^2_R\big]dt\\
&+ \sum_{i=1}^N \mathbb{E}\big[|z_i(T)|_H^2-|z^{(N)}(T)|^2_{H_{\hat{\Gamma}}}
-2( z_0(T))^TH_{\hat{\Gamma}_1}^Tz^{(N)}(T)+|z_0(T)|^2_{\hat{\Gamma}_1^TH\hat{\Gamma}_1}\big].\nonumber
\end{align}
{Let $\{\breve{p}_i, \breve{q}_i^j,i,j=0,1,\cdots,N\}$ be a set of solutions to (\ref{eq59}).
Then, by It\^{o}'s formula, we obtain 
\begin{align*}
	&\sum_{i=1}^N\mathbb{E}\big[\langle\hat{\Gamma}^T_1H (\hat{\Gamma}-I)\breve{x}^{(N)}(T)
	+\hat{\Gamma}_1^TH\hat{\Gamma}_1\breve{x}_0^T(T), z_0(T)\rangle\big]\cr
	&=\sum_{i=1}^N\mathbb{E}[\langle \breve{p}_0(T),z_0(T)\rangle-\langle \breve{p}_0(0),z_0(0)\rangle]\cr
	=&\sum_{i=1}^N\mathbb{E}\!\int_0^T\! \Big\{\big\langle -\big[(A_0+B_0P_0)^T\breve{p}_0+F^T\breve{p}^{(N)}+(C_0+D_0P_0)^T\breve{q}_0^0+\bar{F}^T\breve{q}^{(N)}\cr
	&-\Gamma_1^TQ((I-\Gamma) \breve{x}^{(N)}-\Gamma_1 \breve{x}_0)\big],z_0\big\rangle+\langle \breve{p}_0,(A_0+B_0P_0)z_0+({G}_0+B_0\bar{P})z^{(N)}\rangle\cr
	&+\langle \breve{q}_0^0,(C_0+D_0P_0)z_0+(\bar{G}_0+D_0\bar{P})z^{(N)}\rangle\Big\} dt\cr
	= &\sum_{i=1}^N \mathbb{E}\int_0^T \Big\{\big\langle -[F\breve{p}^{(N)}+\bar{F}\breve{q}^{(N)}-\Gamma_1^TQ((I-\Gamma) \breve{x}^{(N)}-\Gamma_1\breve{x}_0)],z_0\big\rangle\cr
	&+\langle({G}_0+B_0\bar{P})^T\breve{p}_0+(\bar{G}_0+D_0\bar{P})^T\breve{q}^0_0,z_i\big\rangle\Big\} dt,
\end{align*}
and
\begin{align*}
	&\sum_{i=1}^N\mathbb{E}[\langle H\breve{x}_i(T)-H_{\hat{\Gamma}}\breve{x}^{(N)}(T)+(\hat{\Gamma}-I)^TH\hat{\Gamma}_1 \breve{x}_0(T),z_i(T)\rangle]\cr
	= &\sum_{i=1}^N \mathbb{E}\int_0^T \Big\{\big\langle -\big[Q\breve{x}_i-Q_{\Gamma}\breve{x}^{(N)}
	+(\Gamma-I)^TQ\Gamma_1\breve{x}_0+(\bar{G}_0+B_0\bar{P})^T\breve{p}_0+(\bar{G}_0+D_0\bar{P})^T\breve{q}^0_0\big],z_i\big\rangle\cr
	&+\langle F\breve{p}^{(N)}+\bar{F}\breve{q}^{(N)},z_0\rangle+\langle B^T\breve{p}_i+D^T\breve{q}_i^i,u_i\rangle\Big\} dt,
\end{align*}
where
the second equation holds since $\sum_{i=1}^N\mathbb{E}\langle G^T\breve{p}^{(N)},z_i\rangle=\sum_{i=1}^N\mathbb{E}\langle \breve{p}_i,Gz^{(N)}\rangle $
and  $\sum_{i=1}^N\mathbb{E}\langle \bar{G}^T\breve{q}^{(N)},z_i\rangle=\sum_{i=1}^N\mathbb{E}\langle \breve{q}_i^i,\bar{G}z^{(N)}\rangle .$
From the above equations and (\ref{eq60}), 
\begin{align}
	I_1=&\frac{1}{N}\sum_{i=1}^N \mathbb{E}\int_0^T\big[\big\langle Q\breve{x}_i-{Q_{\Gamma}}\breve{x}^{(N)}+(\Gamma-I)^TQ\Gamma \breve{x}_0,z_i\big\rangle+\langle\Gamma_1^T Q(\Gamma-I)\breve{x}^{(N)} +\Gamma_1^TQ\Gamma_1\breve{x}_0,z_0\rangle\cr
	&+ \langle R\breve{u}_i,u_i\rangle\big]dt+ \sum_{i=1}^N \mathbb{E}\big[\big\langle H\breve{x}_i(T)-{H_{\hat{\Gamma}}}\breve{x}^{(N)}(T)+(\hat{\Gamma}-I)^TH\hat{\Gamma} x_0(T),z_i(T)\big\rangle\cr
	&+\big\langle\hat{\Gamma}^T_1H (\hat{\Gamma}-I)\breve{x}^{(N)}(T)
	+\hat{\Gamma}_1^TH\hat{\Gamma}_1\breve{x}_0^T(T), z_0(T)\big\rangle\big]\cr
	=&\frac{1}{N}\sum_{i=1}^N \mathbb{E}\int_0^T[\langle R\breve{u}_i+B^T\breve{p}_i+D^T\breve{q}_i^i,u_i\rangle]dt.
\end{align}
Note that $Q-Q_{\Gamma}=(I-\Gamma)^TQ(I-\Gamma)$ and
$H-H_{\hat{\Gamma}}=(I-\hat{\Gamma})^TH(I-\hat{\Gamma})$. Then, we have
$$\begin{aligned}
	I_2=&\frac{1}{N}\sum_{i=1}^N \mathbb{E}\int_0^T\big[|z_i-z^{(N)}|_Q^2+|z^{(N)}|^2_{Q-Q_{\Gamma}}+2(\Gamma z_0)^TQ(\Gamma-I)z^{(N)}+|\Gamma_1 z_0|^2_Q+|{u}_i|^2_R\big]dt\cr
	&+ \sum_{i=1}^N \mathbb{E}\big[|z_i(T)-z^{(N)}(T)|_H^2+|z^{(N)}(T)|^2_{H-H_{\hat{\Gamma}}}-2 z_0^T(T)H_{\hat{\Gamma}_1}^Tz^{(N)}(T)+|\hat{\Gamma}_1z_0(T)|^2_H\big]\cr
	=&\frac{1}{N}\sum_{i=1}^N \mathbb{E}\int_0^T\big[|z_i-z^{(N)}|_Q^2+|(I-\Gamma)z^{(N)}-\Gamma_1 z_0|^2_{Q}+|{u}_i|^2_R\big]dt
	\cr
	&+ \sum_{i=1}^N \mathbb{E}\big[|z_i(T)-z^{(N)}(T)|_H^2+|(I-\hat{\Gamma})z^{(N)}(T)-\hat{\Gamma}_1z_0(T)|^2_{H}
	\big].
\end{aligned}$$
Since $Q\geq0$, $R>0$, and $H\geq0$,  we obtain
$I_2\geq0$. From (\ref{eq5-3}), $\breve{u}$ is a minimizer to (P1) if and only if
$I_1=0 $, which is equivalent to
$R\breve{u}_i+B^T\breve{p}_i+D^T\breve{q}_i^i=0, \ i=1,\cdots,N.$
Thus, we have the optimality system (\ref{eq59}).
This implies that (\ref{eq59}) admits a solution $(\check{x}_i,\check{p}_i,\check{q}_{i}^{j}, i,j=1,\cdots,N)$.
\hfill{$\Box$}

\emph{Proof of Theorem \ref{thm4.4}.} \emph{(For followers)}. 
By (\ref{eq3b}),
it can be verified that under feedback strategies (\ref{eq3a}), $\mathbb{E}\int_0^T(|\bar{x}_0|^2+|\bar{x}|^2)dt<c$. This  gives
$  \mathbb{E}\int_0^T(|{x}_i|^2+|{x}^{(N)}|^2)dt<c_1.$
 Furthermore, from (\ref{eq3b}) and (\ref{eq75a}) we have
$$\begin{aligned}
	d({x}^{(N)}-\bar{x})=&(A+G+B\hat{K})
	({x}^{(N)}-\bar{x})dt\cr
	&+\frac{1}{N}\sum_{j=1}^N[(C+D\hat{K}){x}_i+\bar{G}{x}^{(N)}+D\hat{K}\bar{x}+(\bar{F}+DK_0){x}_0]dW_j,
\end{aligned}$$
Similar to (\ref{eq41a}), we  have for any $t\in [0,T]$,
\begin{align}\label{eq87}
	\mathbb{E}&|{x}^{(N)}(t)-\bar{x}(t)|^2\leq \big|\bar{\Xi}(t,0)\big|^2\mathbb{E}|{x}^{(N)}(0)-\bar{x}(0)|^2\cr
	&+\frac{1}{N^2}\sum_{i=1}^N\int_0^t c\big|\bar{\Xi}(t,s)\big|
	\max_{1\leq i\leq N}\mathbb{E}\big(|{x}_i|^2+|{x}^{(N)}|^2+|\bar{x}|^2+|{x}_0|^2)ds
	=O(\frac{1}{N}),
\end{align}
where  $\bar{\Xi}(\cdot,s)$ satisfies $\frac{d\bar{\Xi}(t,s)}{dt}=({A}+G+B\hat{K})\bar{\Xi}(t,s),$ $\bar{\Xi}(s,s)=I.$ This further gives $	\mathbb{E}|{x}_0(t)-\bar{x}_0(t)|^2=O(1/N)$,
for any $0\leq t\leq T$.
Note that  $\bar{x}=\mathbb{E}[{x}_i|\mathcal{F}^0]=\mathbb{E}[{x}^{(N)}|\mathcal{F}^0]$ (which follows  from (\ref{eq3b})). Then, we have
\begin{equation}\label{eq86b}
	\mathbb{E}[\bar{x}^T(x^{(N)}-\bar{x})]=\mathbb{E}\big[\bar{x}^T\mathbb{E}[x^{(N)}-\bar{x}|\mathcal{F}^0]\big]=0.
\end{equation}
From (\ref{eq2b}) and (\ref{eq87}), we have
\begin{align}
   \label{eq86}
	J_{\rm soc}^{(N)}(u_0,u)
	=&\frac{1}{N}\sum_{i=1}^N \mathbb{E}\int_0^T\big[|x_i|_Q^2-|x^{(N)}|^2_{Q_{\Gamma}}-2x_0^TQ_{\hat{\Gamma}_1} ^T x^{(N)}+|\Gamma_1x_0|^2_Q+|u_i|^2_R\big]dt\\
	&+ \frac{1}{N}\sum_{i=1}^N \mathbb{E}\big[|x_i(T)|_H^2
	-|x^{(N)}(T)|^2_{H_{\hat{\Gamma}}}-2(H_{\hat{\Gamma}_1} x_0(T))^T\bar{x}(T)+|\Gamma_1x_0(T)|^2_H\big]
\end{align}
\begin{align}
	\leq&\frac{1}{N}\sum_{i=1}^N \mathbb{E}\int_0^T\big[|x_i|_Q^2-|\bar{x}|^2_{Q_{\Gamma}}-2x_0^TQ_{\hat{\Gamma}_1} ^T\bar{x}+|\Gamma_1x_0|^2_Q+|u_i|^2_R\big]dt\cr
	&+\frac{1}{N} \sum_{i=1}^N \mathbb{E}\big[|x_i(T)|_H^2-|\bar{x}(T)|^2_{H_{\hat{\Gamma}}}-2(H_{\hat{\Gamma}_1} x_0(T))^T\bar{x}(T)+|\Gamma_1x_0(T)|^2_H\big]+ \epsilon_1\cr
	\stackrel{\Delta} {=}&\bar{J}^{(N)}_{\rm soc}(u_0,u)+\epsilon_1.\nonumber
\end{align}
We now deform $\bar{J}^{(N)}_{\rm soc}(u_0,u)$ by the method of completing squares. Note that  $\bar{x}=\mathbb{E}[{x}_i|\mathcal{F}^0]$ satisfies
\begin{equation}\label{eq88}
	d\bar{x}=[(A+G)\bar{x}+B\bar{u}+Fx_0]dt,
\end{equation}
where $\bar{u}=\mathbb{E}[{{u}}_i|\mathcal{F}^0]$. Then, it follows that
$$\begin{aligned}
	d(x_i-\bar{x})=&[A(x_i-\bar{x})+B(u_i-\bar{u})+G(x^{(N)}-\bar{x})]dt+(Cx_i+Du_i+\bar{G}x^{(N)}+\bar{F}x_0)dW_i.
\end{aligned}$$
From (\ref{eq86b}), applying It\^{o}'s formula to $|x_i-\bar{x}|^2_{M}$, we obtain
\begin{align}\label{eq89}
	&\mathbb{E}\big[|x_i(T)-\bar{x}(T)|^2_{H}-|x_i(0)-\bar{x}(0)|^2_{M(0)}\big] \\
	=&\mathbb{E}\int_0^T\Big\{ (x_i-\bar{x})^T(\dot{M}+A^TM+MA+C^TMC)(x_i-\bar{x})+(u_i-\bar{u})^TD^TMD(u_i-\bar{u})\cr
	&+2(u_i-\bar{u})^T(B^TM+D^TMC)(x_i-\bar{x})+\bar{u}^TD^TMD\bar{u}+x_0^T\bar{F}^TM\bar{F}x_0\cr
	&+\bar{x}^T(C+G)^TM[(C+\bar{G})\bar{x}+2\bar{F}x_0]+2\bar{u}^TD^TM[(C+\bar{G})\bar{x}+\bar{F}x_0]\cr
	&+2(x^{(N)}-\bar{x})^T[(\bar{G}^TMC+G^TM)(x_i-\bar{x})+\bar{G}^TMD(u_i-\bar{u})]\Big\}dt.\nonumber
\end{align}
It follows by (\ref{eq88}) that
\begin{align}
	&\mathbb{E}\big[\bar{x}^T(T)(H-H_{\hat{\Gamma}})\bar{x}(T)-\bar{x}^T(0)(M(0)+\bar{M}(0))\bar{x}(0)\big]\\
	=&\mathbb{E}\int_0^T \big\{ \bar{x}^T[\dot{M}+\dot{\bar{M}}+(A+G)^T(M+\bar{M})+(M+\bar{M})(A+G)]\bar{x}\cr
	&+2\bar{x}^T(M+\bar{M})B\bar{u}+2\bar{x}^T(M+\bar{M})Fx_0\big\}dt.\nonumber
\end{align}
By (\ref{eq3b}) and It\^{o}'s formula,
\begin{align}
	&\mathbb{E}\big[{x}_0^T(T)\hat{\Gamma}_1^TH\hat{\Gamma}_1{x}_0(T)-{x}_0^T(0)\Lambda^0(0){x}_0(0)\big]\\
	=&\mathbb{E}\int_0^T \big\{{x}_0^T[\dot{\Lambda}^0+(A_0+B_0P_0)^T\Lambda^0+\Lambda^0(A_0+B_0P_0)+(C_0+D_0P_0)^T\Lambda^0(C_0+D_0P_0)]{x}_0\cr
	&+2{x}_0^T[\Lambda^0(G_0+ B_0\bar{P})+(C_0+D_0P_0)^T\Lambda^0(\bar{G}_0+D_0\bar{P})\big]\bar{x}+\bar{x}^T(\bar{G}_0+D_0\bar{P})^T\Lambda^0 (\bar{G}_0+D_0\bar{P})\bar{x}\big\}dt.\nonumber
\end{align}
Applying It\^{o}'s formula to $x_0^T\bar{\Lambda}\bar{x}$ and $\bar{x}^TM^0{x}_0$, we have
\begin{align}
	&\mathbb{E}\big[-{x}_0^T(T)H_{\hat{\Gamma}_1}^T\bar{x}(T)-{x}_0^T(0)\bar{\Lambda}(0)\bar{x}(0)\big]\\
	=&\mathbb{E}\int_0^T \big\{{x}_0^T[\dot{\bar{\Lambda}}+\bar{\Lambda}(A+G)+(A_0+B_0P_0)^T\bar{\Lambda}]\bar{x}+{x}_0^T\bar{\Lambda}(B\bar{u}+F{x}_0) +\bar{x}^T(G_0+B_0\bar{P})^T\bar{\Lambda}\bar{x}\big\}dt,\nonumber
\end{align}
and
\begin{align}\label{eq93}
	&\mathbb{E}\big[-\bar{x}^T(T)H_{\hat{\Gamma}_1}{x}_0(T)-\bar{x}^T(0)M^0(0){x}_0(0)\big]\\
	=&\mathbb{E}\int_0^T \big\{\bar{x}^T[\dot{M}^0+(A+G)^TM^0+M^0(A_0+B_0P_0)]\bar{x}+(B\bar{u}+F{x}_0)^TM^0{x}_0 +\bar{x}^TM^0(G_0+B_0\bar{P})\bar{x}\big\}dt.\nonumber
\end{align}
From (\ref{eq73}), 
(\ref{eq89})-(\ref{eq93}), one can obtain
\begin{align*}
	&\bar{J}^{(N)}_{\rm soc}(u_0,u)\cr
	=&\frac{1}{N}\sum_{i=1}^N \mathbb{E}\int_0^T\big[|x_i-\bar{x}|_Q^2+|\bar{x}|^2_{Q-Q_{\Gamma}}+2[(\Gamma-I)^TQ\Gamma_1 x_0]^T\bar{x}+|\Gamma_1x_0|^2_Q+|u_i-\bar{u}|^2_R+|\bar{u}|^2_R\big]dt\cr
	&+\frac{1}{N}\sum_{i=1}^N \mathbb{E}\big[|x_i(T)-\bar{x}(T)|_H^2+|\bar{x}(T)|^2_{H-H_{\hat{\Gamma}}}+2[(\hat{\Gamma}-I)^TH\hat{\Gamma}_1 x_0(T)]^T\bar{x}(T)+|\Gamma_1x_0(T)|^2_H\big]\cr
	=&\frac{1}{N}\sum_{i=1}^N \mathbb{E}\big[|x_i(0)-\bar{x}(0)|^2_{M(0)}+|\bar{x}(0)|^2_{M(0)
		+\bar{M}(0)}+2x_0^T(0)\bar{\Lambda}(0)x^{(N)}(0)+|x_0(0)|^2_{\Lambda_0(0)}\big]\cr
	&+\frac{1}{N}\sum_{i=1}^N \mathbb{E}\int_0^T\Big\{(x_i-\bar{x})^T\Psi^T\Upsilon^{-1}\Psi(x_i-\bar{x})+(u_i-\bar{u})^T\Upsilon(u_i-\bar{u})+2(u_i-\bar{u})^T\Psi(x_i-\bar{x})\cr
	&+\bar{u}^T\Upsilon\bar{u}+\bar{x}^T(\Psi+\bar{\Psi})^T\Upsilon^{-1}(\Psi+\bar{\Psi})\bar{x}
	+2\bar{u}^T[(\Psi+\bar{\Psi})\bar{x}+\Psi^0x_0]+(\Psi^0 x_0)^{T}\Upsilon^{-1}\Psi^0 x_0\cr
	&+2\bar{x}^T(\Psi+\bar{\Psi})^T\Upsilon^{-1}\Psi^0x_0+2(x^{(N)}-\bar{x})^T[(\bar{G}^TMC+G^TM)(x_i-\bar{x})+\bar{G}^TMD(u_i-\bar{u})]\Big\}dt\cr
	=&\frac{1}{N}\sum_{i=1}^N \mathbb{E}\big[|\xi_i|^2_{M(0)}+|\bar{\xi}|^2_{\bar{M}(0)}+2\xi_0^T\bar{\Lambda}(0)\xi_i+|\xi_0|^2_{\Lambda_0(0)}\big]\cr
	&+\frac{1}{N}\sum_{i=1}^N \mathbb{E}\int_0^T\Big\{|u_i-\bar{u}+\Upsilon^{-1}\Psi(x_i-\bar{x})|^2_{\Upsilon}+|\bar{u}+\Upsilon^{-1}[(\Psi+\bar{\Psi})\bar{x}+\Psi^0 x_0]|^2_{\Upsilon}\cr
	&+2(x^{(N)}-\bar{x})^T[\bar{G}^TMC+G^TM](x_i-\bar{x})+\bar{G}^TMD(u_i-\bar{u})]\Big\}dt\cr
	\geq&\frac{1}{N}\sum_{i=1}^N \mathbb{E}\big[|\xi_i|^2_{M(0)}+|\bar{\xi}|^2_{\bar{M}(0)}+2\xi_0^T\bar{\Lambda}(0)\xi_i+|\xi_0|^2_{\Lambda_0(0)}\big]\cr
	&+\frac{1}{N}\sum_{i=1}^N \mathbb{E}\int_0^T2(x^{(N)}-\bar{x})^T[(\bar{G}^TMC+G^TM)(x_i-\bar{x})+\bar{G}^TMD(u_i-\bar{u})]dt.
\end{align*}
Note that  $\hat{u}_i=-\Upsilon^{-1}(\Psi x_i+\bar{\Psi}\bar{x}+\Psi^0x_0)$. From (\ref{eq87}) and (\ref{eq86}), we have
${J}^{(N)}_{\rm soc}(\hat{u}_0,\hat{u})\leq J_{\rm soc}^{(N)}(\hat{u}_0,u)+\epsilon_1,$
where $\epsilon_1=O(1/\sqrt{N})$.

\emph{(For the leader).} 
From (\ref{eq2a}), we have
\begin{align}\label{eq36a}
	J_0(\hat{u}_0,\hat{u}(\hat{u}_0)) \leq& \bar{J}_0(\hat{u}_0,\hat{u}(\hat{u}_0))+\mathbb{E}\int_{0}^{T}\Big[2\big(|x_0(t)- \Gamma_0\bar{x}(t)|^2|Q_0\Gamma_0(\hat{x}^{(N)}(t)-\bar{x}(t))|^2\big)^{1/2}\\
	&+|\Gamma_0(\hat{x}^{(N)}(t)-\bar{x}(t))|_{Q_0}^{2}\Big]dt+|\hat{\Gamma}_0(\hat{x}^{(N)}(T)-\bar{x}(T))|_{H_0}^{2}\Big]\cr
	&+2\mathbb{E}\Big[\big(|x_0(T)- \hat{\Gamma}_0\bar{x}(T)|^2|H_0\hat{\Gamma}_0(\hat{x}^{(N)}(T)-\bar{x}(T))|^2\big)^{1/2}\cr
	\leq &\bar{J}_0(\hat{u}_0,\hat{u}(\hat{u}_0))+O(1/\sqrt{N}).\nonumber
\end{align}
By It\^{o}'s formula, one can obtain
\begin{align}\label{eq32b}
	&\mathbb{E}[x_0^T(T)H_0x_0(T)]-  \mathbb{E}[x_0^T(0)\Theta_1(0)x_0(0)]\\
	=&\mathbb{E} \int_0^T\big[x_0^T(\dot{ \Theta}_1+A_0^T\Theta_1+\Theta_1A_0+C_0^T\Theta_1C_0)x_0+2u_0^T(B^T_0\Theta_1+D_0^T \Theta_1C_0)x_0\big]dt,\label{eq33} \cr
	&\mathbb{E}[\bar{x}^T(T)\hat{\Gamma}_0^TH_0\hat{\Gamma}_0\bar{x}(T)]-  \mathbb{E}[\bar{x}^T(0)\Theta_2(0)\bar{x}(0)]
	=\mathbb{E} \int_0^T\big[\bar{x}^T(\dot{\Theta}_2+\hat{A}^T\Theta_2+\Theta_2\hat{A})\bar{x}+2x_0^T\hat{F}^T\Theta_2\bar{x}\big]dt,\nonumber
\end{align}
and
\begin{align}\label{eq34}
	&\mathbb{E}[\bar{x}^T(T)(-\hat{\Gamma}_0^TH_0){x}_0(T)]-  \mathbb{E}[\bar{x}^T(0)\Theta_3(0){x}_0(0)]\\
	=&\mathbb{E} \int_0^T\big[\bar{x}^T(\dot{\Theta}_3+\hat{A}^T\Theta_3+\Theta_3{A}_0){x}_0+\bar{x}^T\Theta_3B_0u_0+x_0^T\hat{F}^T\Theta_3{x}_0\big]dt.\nonumber
\end{align}
It follows from (\ref{eq32b})-(\ref{eq34}) that
\begin{align}
	\bar{J}_0(u_0,{u}(u_0))=& \mathbb{E}[x_0^T(0)\Theta_1(0)x_0(0)+\bar{x}^T(0)\Theta_2(0)\bar{x}(0)+\bar{x}^T(0)\Theta_3(0){x}_0(0)]\\
	&+\mathbb{E}\int_0^T\Big[x^T_0(B_0^T\Theta_1+D_0^T\Theta_1C_0)^T\Xi^{-1}(B_0^T\Theta_1+D_0^T\Theta_1C_0)x_0\cr
	&+\bar{x}^T\Theta_3B_0\Xi^{-1}B^T_0\Theta_3\bar{x}+2\bar{x}^T\Theta_3B_0\Xi^{-1}(B_0^T\Theta_1+D_0^T\Theta_1C_0)x_0 \cr
	&+2 u_0^T[(B_0^T\Theta_1+D_0^T\Theta_1C_0)x_0+B_0^T\Theta_3\bar{x}]+u_0^T\Xi u_0 \Big]dt\cr
	=& \mathbb{E}[\xi_0^T\Theta_1(0)\xi_0+\bar{\xi}^T\Theta_2(0)\bar{\xi}+\bar{\xi}^T\Theta_3(0)\xi_0]
	+\mathbb{E}\int_0^T\Big[\big|u_0\cr
	&+\Xi^{-1}(B_0^T\Theta_1+D_0^T\Theta_1C_0)x_0+\Xi^{-1}B_0^T\Theta_3\bar{x}\big|^2_{\Xi}\Big]dt\cr
	\geq &\mathbb{E}[\xi_0^T\Theta_1(0)\xi_0+\bar{\xi}^T\Theta_2(0)\bar{\xi}+\bar{\xi}^T\Theta_3(0)\xi_0]
	=\bar{J}_0(\hat{u}_0,\hat{u}(\hat{u}_0)).\nonumber
\end{align}
This together with (\ref{eq36a}) leads to $   J_0(\hat{u}_0,\hat{u}(\hat{u}_0))\leq  \bar{J}_0(u_0,u^*({u}_0))+O(1/\sqrt{N}).$
The reminder of the proof is similar to that of Theorem \ref{thm4.2}. \hfill{$\Box$}

\end{document}